\documentclass[final,reqno,onefignum,onetabnum]{siamltex}
\usepackage{enumerate}
\usepackage{amsmath,amssymb,graphicx,bbm}
\usepackage{verbatim,rotating}
\usepackage{mathrsfs,fancyhdr,xcolor}

\newtheorem{Theorem}{\indent Theorem}[section]
\usepackage[ruled,vlined]{algorithm2e}
\usepackage{subfig}
\newtheorem{Remark}{\indent Remark}[section]

\DeclareMathOperator*{\argmin}{argmin}

\def\R{\mathbb{R}}

\title
{
{Data-driven polynomial chaos expansions: a weighted least-square approximation}\thanks{The first author is
supported by the NSF of China (No.11671265) and Program for Outstanding Academic leaders in Shanghai (No.151503100). The last author is partially supported by the NSF of China (under grant numbers 11688101, 91630312, 91630203, 11571351, and 11731006), the science challenge project (No. TZ2018001), NCMIS, and the youth innovation promotion association (CAS).}
}

\author
{
Ling Guo\thanks{Department of Mathematics, Shanghai Normal University, Shanghai, China. Email: lguo@shnu.edu.cn.}
\and
Yongle Liu\thanks{Department of Mathematics, Southern University of Science and Technology, Shenzhen, China. Email: 11749318@mail.sustc.edu.cn.}
\and Tao Zhou\thanks{LSEC, Institute of Computational Mathematics and Scientific/Engineering
Computing, Academy of Mathematics and Systems Science, Chinese Academy of Sciences, Beijing, China. Email: tzhou@lsec.cc.ac.cn.}
}

\begin{document}
\maketitle

\graphicspath{{figure/}}

\begin{abstract}
In this work, we combine the idea of  data-driven polynomial chaos expansions with the weighted least-square approach to solve uncertainty quantification (UQ) problems. The idea of data-driven polynomial chaos is to use statistical moments of the input random variables to develop an arbitrary polynomial chaos expansion, and then use such data-driven bases to perform UQ computations. Here we adopt the bases construction procedure by following \cite{Ahlfeld_2016SAMBA}, where the bases are computed by using matrix operations on the Hankel matrix of moments. Different from previous works, in the postprocessing part, we propose a weighted least-squares approach to solve UQ problems. This approach includes a sampling strategy and a least-squares solver. The main features of our approach are two folds: On one hand, our sampling strategy is independent of the random input. More precisely, we propose to sampling with the equilibrium measure, and this measure is also independent of the data-driven bases. Thus, this procedure can be done in prior (or in a off-line manner). On the other hand, we propose to solve a Christoffel function weighted least-square problem, and this strategy is quasi-linearly stable -- the required number of PDE solvers depends linearly (up to a logarithmic factor) on the number of (data-driven) bases. This new approach is thus promising in dealing with a class of problems with epistemic uncertainties. Several numerical tests are presented to show the effectiveness of our approach.
\end{abstract}



\pagestyle{myheadings}
\thispagestyle{plain}

\section{Introduction}
Uncertainty Quantification (UQ) has been a hot topic recently. The aim of UQ is to quantify the impact of the stochastic inputs to the stochastic response, and thus a fundamental problem of UQ is to approximate a potentially high dimensional parametric function $f(\xi_1,\xi_2,...,\xi_d): \mathbb{R}^d\rightarrow\mathbb{R}, \, d\geq1.$  One popular way to perform UQ analysis is to assume that  the distributions of the input parameters $\{\xi_k\}_{k=1}^d$ are known in prior, and this is also well known as aleatory-type uncertainty model. Among others, the generalized Polynomial Chaos (gPC) \cite{Xiu_2002Wiener} based on the Wiener-Askey formula, which is an extension of the original work by Wiener \cite{Wiener_1938homochaos}, is a popular approach for aleatory-type uncertainty analysis.  The idea is to approximate the parametric function $f$ with polynomial bases that are orthogonal with respect to the input density of the parameters. The unknown expansion coefficients can then be computed by performing for example the Galerkin projection into a finite polynomial space. Notice that in general, one needs to solve a coupled Galerkin system that is much more complicated than the original model -- the so called intrusive approach. Another popular approach, termed stochastic collocation, has gained much attention due to its efficiency and its non-intrusive property. The idea of stochastic collocation is to use efficient sample solutions to construct global polynomial approximations. For recent developments of stochastic collocation methods, one can refer to \cite{Narayan_2015SCUMM,Nobile_2008sparse,Tang_2010ConSc,Tang_2015redeUQ} and references therein.

In recent years, there is a growing demand to extend the gPC method to more general input distributions (beyond the Wiener-Askey formula). One of the first attempts is the multi-element generalized polynomial chaos (ME-gPC) \cite{Prempraneerach_2010Uq,Wan_2008MEgPC} where the random space is divided into small elements and local polynomial expansion is constructed via the Gram-Schmidt procedure. While the approach gains advantages when dealing with discontinuity of model responses, the computational complexity can be very dramatic. A multi-element probabilistic collocation method was also developed along this direction \cite{Zheng_2015MEPC}. Global polynomial expansions for arbitrary distributions have also been investigated based on Gram-Schmidt orthogonalisation \cite{Witteveen_07uncerapc,Witteveen_2006MauGS}. However, such approaches still rely on the availability of the input density function.

More recently,  Oladyshkin and Nowak \cite{Oladyskin_2012Datadriven} propose a moment match method to deal with arbitrary distributions (termed aPC), and the approach is promising when one has incomplete input information, such as the situation when only sample locations are given.  The idea in \cite{Oladyskin_2012Datadriven} is to set up the moment match equations, and then solve the unknown polynomial coefficients.  The aPC offers a possibility to propagate only the given information without making assumptions. As showed by Oladyshkin and Nowak in \cite{Oladyshkin_2011concept} that only moments are propagated in all PC approaches, thus the aPC offers the most reliable results with limited input data. Although the aPC construction approach are straightforward to implement, it is well known that the coefficient matrix of the moment equation maybe ill conditioned when the polynomial order is large. Recently, a promising alternative way to calculate the aPC was proposed in \cite{Ahlfeld_2016SAMBA}, where the authors proposed an algorithm in which all the required quantities are calculated directly using only matrix operations performed on the Hankel matrix of moments. Then,  a sparse grid approach based on the Smolyak¡¯s algorithm was proposed in \cite{Ahlfeld_2016SAMBA} where the collocation points are generated by the constructed bases, yet again by using matrix operations.

Unlike the traditional gPC methods, where one perform UQ computations directly based on well known polynomial bases choosing according to the Wiener-Askey formula, the aPC approach can normally be divided into the following two steps:
\begin{itemize}
\item Bases construction. One uses the input information (moments, samples locations, ect.)  to construct the so called arbitrary polynomial bases (data-driven bases). Notice that this procedure is somehow model-independent and only input information is used.

 \item UQ computations. One adopts the data-driven bases to perform UQ computations. This procedure is obviously model-dependent, and one could consider a stochastic Galerkin approach, or a sparse grid stochastic collocation approach as in \cite{Oladyskin_2012Datadriven} (where collocation points are generated using the arbitrary polynomial bases).
\end{itemize}

In this work, the only information we needed are some sample locations (The density of the input is unknown). We shall then adopt the aPC construction procedure in \cite{Ahlfeld_2016SAMBA}.  However, in the second (postprocessing) step, we propose a weighted least-squares approach to obtain the aPC expansion coefficients. This approach includes a sampling strategy and a least-squares solver. We propose to sampling with the equilibrium measure which is independent of the data driven bases (or the input information). Thus, this procedure can be done in prior (or in a off-line manner). Then we propose to solve a Christoffel function weighted least-squares problem, and in many cases of interests this approach is  linearly stable -- the number of samples (the number of PDE solvers) depends linearly on the number of (data-driven) bases. We shall present theoretical motivations and several numerical tests to support our statements.

The rest of this paper is organized as follows. In Section 2, we introduce the traditional gPC approach. The construction procedure of  data-driven polynomial bases is introduced in Section 3. In Section 4, we present a weighted least-squares approach to perform UQ computations.  Numerical experiments are then shown in Section 5 to indicate the applicable and effectiveness of our approach. Finally, we give some concluding remarks in Section 6.

\section{Generalized polynomial chaos}

In parametric uncertainty quantification studies, the main goal is to trace the effect of the random inputs, here denoted by $\xi=(\xi_1,\xi_2,\ldots,\xi_d)$ through the model and to quantify their effect on the model output (prediction) $f(\xi): \mathbb{R}^d\rightarrow \mathbb{R}$. This is frequently done via the generalized polynomial chaos expansions. Concretely, we assume that the components of the random input $\xi=(\xi_1,\xi_2,\ldots,\xi_d)$ are mutually
independent, and for each $\xi_i$ in $\Gamma_i \subset \R$ it admits a marginal probability density $\rho_i$. Then the joint density function for $\xi$ yields $\rho(\xi)= \prod_{i=1}^d \rho_i(\xi_i): \Gamma\rightarrow \mathbb{R}^+$ with $\Gamma:= \prod_{i=1}^d\Gamma_i \subset \mathbb{R}^d.$  The gPC approach seeks to construct a polynomial approximation of
 $f(\xi)$ as follows:
\begin{equation}\label{eq:multi-expan}
f(\xi)\approx\sum_{\boldsymbol{\alpha}\in\Lambda}c_{\boldsymbol{\alpha}}\Phi_{\boldsymbol{\alpha}}(\xi),
\end{equation}
where $\boldsymbol{\alpha}=\{{\alpha}_1,{\alpha}_2,\ldots,{\alpha}_d\}$ is a multi-index and $\Lambda$ is a finite multi-index set. And $\Phi_{\boldsymbol{\alpha}}$ is the multivariate orthogonal polynomials that are orthogonal with respect to the density $\rho(\xi),$ i.e.,
\begin{equation}
\int_\Gamma \!\! \rho(\xi) \Phi_{\boldsymbol{\alpha}}(\xi) \Phi_{\boldsymbol{\beta}}(\xi)d\xi = \delta_{\boldsymbol{\alpha},\boldsymbol{\beta}}, \quad \boldsymbol{\alpha},\boldsymbol{\beta}\in\Lambda.
\end{equation}
Notice the polynomials are defined as tensor-products of the univariate orthogonal polynomials in each direction, i.e.,
\begin{equation*}\label{eq:d-variate-polynomial}
\Phi_{\boldsymbol\alpha}=\prod_{i=1}^d\phi_{{\alpha}_i}^i(\xi_i) \quad  \textmd{with} \quad
  \int_{\Gamma_i}\!\! \phi^i_{{\alpha}_k}(\xi_i) \phi^i_{{\alpha}_l}(\xi_i) \rho_i(\xi_i) d\xi_i = \delta_{k,l}.
\end{equation*}
In this work we focus on the total degree polynomial space that is defined as
\begin{align}\label{eq:TD}
P(\Lambda)=\mathrm{span}\left\{ \Phi_{\boldsymbol\alpha} \;\; \big| \;\;  \boldsymbol\alpha \in \Lambda_k^{\textrm{TD}}, \quad \textmd{with} \quad \Lambda_k^{\textrm{TD}}:=\left\{\boldsymbol\alpha \;\; \big| \;\;  |\boldsymbol\alpha|_1=\sum_{i=1}^d \alpha_i \leq k \right\}\right\}.
\end{align}
It is usually more convenient use the single index instead of the multi-index, and to this end, one can place an order on the multi-indices, i.e.,
\begin{align}\label{eq:index-ordering}
  \left\{ \boldsymbol{\alpha} \;\; | \;\; \boldsymbol{\alpha} \in \Lambda \right\} \longleftrightarrow \left\{ 1, \ldots, N \right\}.
\end{align}
Thus we have
\begin{equation}\label{eq:multi-index-single}
\{\Phi_{\boldsymbol{\alpha}}(\xi)\}_{\boldsymbol{\alpha}\in\Lambda}\Leftrightarrow \{\Phi_j(\xi)\}_{j=1}^{N=\dim(P(\Lambda))}.
\end{equation}
Hereafter, for simplicity, we will use the single index $\{j=1,2,\ldots,N\}$. Therefore, the gPC approximation (\ref{eq:multi-expan}) can be written as
\begin{equation}\label{eq:finite-N-expan}
f(\xi)\approx f_N(\xi)=\sum_{j=1}^Nc_j\Phi_j(\xi).
\end{equation}
The main purpose now is to estimate the coefficients $\{c_j\}_{j=1}^N$ in an efficient way. Many numerical techniques on how to obtain the polynomial coefficients in UQ problems have been developed in recent years, such as the intrusive stochastic Galerkin methods \cite{Tang_Zhou_JSC,MNZ_SISC,Xiu_2002Wiener} and the non-intrusive collocation methods \cite{Eldred_2009nonintrusive,Nobile_2008sparse,Zhou_2015wlsrq,XiuH,Zhou_2014multilsweil,Jakeman_2016generalizedsample,Gao_2014lsnodels,Guo_2018}.

\section{Data-driven polynomial chaos: a moment-based approach }
The gPC methods discussed above assume an exact knowledge of the involved probability density functions. However, the distribution information of the random input is very limited in many engineering applications, and sometimes the only information available is sample locations. This is also unknown as epistatic uncertainty. To deal with these situations, more general types of polynomial chaos expansions have been investigated in the past few years, see e.g. \cite{Witteveen_07uncerapc,Witteveen_2006MauGS,Ernst_2012gPCconvergen,
soize_2004physicalrandom,Oladyskin_2012Datadriven,Ahlfeld_2016SAMBA}. Here we shall review the idea of arbitrary polynomial chaos (aPC for short) approach developed in \cite{Oladyskin_2012Datadriven,Ahlfeld_2016SAMBA}. Such an approach can handle the situation when one only has sample locations (or only moments information is available for the random input).

\subsection{Moment match approaches}
In this section, we shall review the basic idea in \cite{Oladyskin_2012Datadriven,Ahlfeld_2016SAMBA}. We suppose that we are given moments information for the random input (while the associated distributions are unknown). Notice that this approach provides the possibility to propagate continuous or discrete probability density functions and also histograms (data sets) as long as their moments exist and the determinant of the moment matrix is strictly positive (see details below). The aim is to construct a set of polynomials bases $\{\Phi_j\}$ that admit a good approximation for the underlining parametric problem. This will be done by using the moment match methods. We first present the idea in the one dimensional setting.

Suppose that the density function for a continuous random variable $\eta \in I$ is $\rho(\eta)$, then the $k$-th raw moment $\mu_k$ is defined by
\begin{equation}
\mu_{k}=\int_{I}\eta^{k}\rho(\eta)d\eta, \ \ \ k=0,1,\ldots .
\end{equation}
Similarly, if the random variable $\eta$ is of discrete-type $\eta\in\widehat{I}$ then its $k$-th moment is defined as
\begin{equation}
\mu_{k}=\sum_{\eta\in\widehat{I}}\eta^{k}\rho(\eta), \ \ \ k=0,1,\ldots.
\end{equation}
Finally, if a random variables is only presented as a set of $M$ samples locations $\{\eta_1,\eta_2,\ldots,\eta_M\}$ (The setting in this work), the $k$-th moment $\mu_k$ can be calculated approximately by
\begin{equation}\label{eq:samples}
\mu_{k}=\frac{1}{M}\sum_{m=1}^M\eta_m^k, \ \ \ k=0,1,\ldots.
\end{equation}
Suppose we know the moments of $\eta$ up to the index $2K,$ then we can consider to construct a set of orthogonal polynomial bases    $\{\phi_k(\eta)\}_{k=0}^{K}$ with the general form
$$\phi_k(\eta)=\sum_{j=0}^k\beta_j\eta^j, \quad k=0, ..., K.$$
By matching the moments information, we obtain
\begin{gather}\label{gsshiqi}
\begin{bmatrix}
\mu_{0}&\mu_{1}&\cdots&\mu_{k}\\
\mu_{1}&\mu_{2}&\cdots&\mu_{k+1}\\
\vdots&\vdots&\vdots&\vdots\\
\mu_{k-1}&\mu_{k}&\cdots&\mu_{2k-1}\\
0&0&\cdots&1
\end{bmatrix}
\begin{bmatrix}
\beta_{0}\\
\beta_{1}\\
\vdots\\
\beta_{k-1}\\
\beta_{k}
\end{bmatrix}=\begin{bmatrix}
0\\
0\\
\vdots\\
0\\
1
\end{bmatrix}.
\end{gather}
Thus one can obtain the polynomial coefficients by inverting the above Vandermonde matrix. However, this matrix may become very ill-conditioned when $k$ becomes large. An alternative approach by considering matrix operations on the Hankel matrix of moments was proposed in \cite{Ahlfeld_2016SAMBA}. To introduce the idea, we first define the Hankel matrix of moments as
\begin{gather}\label{eq:Hankelmatrix}
\mathbf{H}=
\begin{bmatrix}
\mu_{0}&\mu_{1}&\cdots&\mu_{k}\\
\mu_{1}&\mu_{2}&\cdots&\mu_{k+1}\\
\vdots&\vdots&\vdots&\vdots\\
\mu_{k}&\mu_{k+1}&\cdots&\mu_{2k}
\end{bmatrix}.
\end{gather}
If the moments are given by samples (\ref{eq:samples}), we require that the set of $M$ samples is determinate in the Hamburger sense, meaning that all the corresponding quadratic forms are strictly positive, that is $\text{det}(\mathbf{H})>0$. Given the above Hankel matrix of moments, we first perform the Cholesky decomposition to obtain $\mathbf{H}=\mathbf{R}^{\top}\mathbf{R}$ with
\begin{gather}\label{eq:Cholmatrix}
\mathbf{R}=
\begin{bmatrix}
r_{11}&r_{12}&\cdots&r_{1,k+1}\\
      &r_{22}&\cdots&r_{2,k+1}\\
      &  &\ddots&\vdots\\
      &   &     &r_{k+1,k+1}
\end{bmatrix}.
\end{gather}
Then, the Mysovskih theorem \cite{Mysovskikh_1968cubature} states that the entries of the matrix $\mathbf{R}$ can form an orthogonal system of polynomials. Moreover, explicit analytic formulas to obtain the polynomial coefficients are available \cite{Golub_1968Gqr}:
\begin{align}\label{eq:three-termrec}
\eta\phi_{j-1}(\eta)=b_{j-1}\phi_{j-2}(\eta)+a_j\phi_{j-1}(\eta)+b_j\phi_j(\eta), \quad j =1, ... k.
\end{align}
Here $a_j$ and $b_j$ can be computed by the components of matrix $\mathbb{R}$:
\begin{align}
a_j=\frac{r_{j,j+1}}{r_{j,j}}-\frac{r_{j-1,j}}{r_{j-1,j-1}}, \quad b_{j}=\frac{r_{j+1,j+1}}{r_{j,j}},
\end{align}
where $r_{0,0}=1$ and $r_{0,1}=0$.

\begin{Remark}
In the above discussions, we have only presented the one dimensional case. For high dimensional cases, one can simply perform the similar procedure as above, and then obtain the multi-variate bases by using the tensor-product rule. Given such data-driven (or moment driven) polynomial bases, one can then perform UQ computations for the underline models. For example, a sparse grid method  was proposed in \cite{Ahlfeld_2016SAMBA}, where the stochastic collocation points are generated again by using matrix operations based on the data-driven bases discussed above.
\end{Remark}

\begin{Remark}
We remark again that the above aPC approach provides the possibility to propagate continuous or discrete probability density functions and also data sets as long as their moments exist and the determinant of the moment matrix is strictly positive. The expansion bases here are fully data-driven, and we do not require any distribution information. For cases with limited data, such an approach can avoid bias and fitting errors caused by wrong assumptions.
\end{Remark}

\subsection{Some theoretical discussions}

We have reviewed the moment match approach for constructing data-driven bases for UQ studies, by requiring that the moment problem is uniquely solvable. Following closely \cite{Ernst_2012gPCconvergen}, we now provide with some mild conditions that can guarantee such an requirement. Our basic
assumptions are as following:
\begin{itemize}
\item Assumption 1: we assume that each basic random variable $\eta$ possesses finite moments of all orders.
\item Assumption 2: the associated distribution functions $F_{\eta}(x) := P(\eta \leq x)$ of the basic random variables are continuous.
\end{itemize}
Notice that such assumptions are just for theoretical analysis, the approach above can still be used even if the probability density functions are of discrete type as long as their moments exist and the determinant of the moment matrix is strictly positive. In other words, the following theorem only works when the input random variables satisfy the above two assumptions. For more general settings, the relevant theoretical foundation is still open.

\begin{Theorem}[\cite{Ernst_2012gPCconvergen}]
If one of the following conditions is valid, then the moment problem is uniquely solvable and therefore the set of polynomials (that constructed by the moment match approach) in the random variable $\eta$ is dense in the space $L^2(\Omega,\sigma(\eta), P)$, where $\Omega$ is the abstract set of elementary events, $\sigma(\eta)$ is a $\sigma$-algebra of subsets
of $\Omega$ and $P$ is  a probability measure on $\sigma(\eta)$.
\begin{enumerate}
  \item The distribution $F_{\eta}$ has compact support, i.e., there exists a compact interval $[a, b], a, b\in\mathbb{R}$, such that $P(\eta\in[a, b]) = 1$.
  \item The moment sequence $\{\mu_k\}_{k\in\mathbb{N}_0}$ of the distribution satisfies
  \begin{equation*}
  \lim_{k\rightarrow\infty}\inf\frac{\sqrt[2k]{\mu_{2k}}}{2k}<\infty.
  \end{equation*}
 \item The random variable is exponential integral,i.e., there holds
 \begin{equation*}
 \langle\exp(a|\eta|)=\int_{\mathbb{R}}\exp{a|x|}F_{\eta}(dx)\rangle<\infty.
 \end{equation*}
 for a strictly positive number $a$. An equivalent condition is the existence of a finite moment-generating function in a neighbourhood of the origin.
  \item (Carleman's condition) The moment sequence $\{\mu_{k}\}_{k\in\mathbb{N}_0}$ of the distribution satisfies
  \begin{equation*}
  \sum_{k=0}^{\infty}\frac{1}{\sqrt[2k]{\mu_{2k}}}=\infty.
  \end{equation*}
 \item (Lin's condition) If the distribution has a symmetric, differentiable and strictly positive density $f_{\eta}$ and for a real number $x_0>0$ there holds
     \begin{equation*}
     \int_{-\infty}^{\infty}\frac{-\log f_{\eta}(x)}{1+x^2}dx=\infty \quad and \quad \frac{-xf_{\eta}^{\prime}(x)}{f_{\eta}(x)}\nearrow \infty (x\rightarrow\infty,x\geq x_0)
     \end{equation*}
\end{enumerate}
\end{Theorem}
The theorem above states that the orthogonal polynomials form a complete bases in $L^2(\Omega,\sigma(\eta), P)$ and thus one can expect a good approximation property using such bases.

\section{Weighted least-squares for postprocessing}

As mentioned above, once we have the data-driven bases, one can perform UQ studies based on such bases. A sparse grid method was proposed in \cite{Ahlfeld_2016SAMBA}, where the collocation points are generated based on the data-driven bases. In this section, we shall propose to use the least-squares approach to do postprocessing computations. Our approach admits many advantages. First of all, we simply sampling with a known measure (the equilibrium measure) to generate collocation points, and the sampling strategy is very cheap and no matrix operations are needed compared to the spares grid approach in \cite{Ahlfeld_2016SAMBA}. Secondly, our sampling strategy is independent of the data-driven bases, and thus this procedure can be done in advance. Finally, our least-squares solver is linear stable in many cases of interests. Details of our approach are presented in the following subsections.

\subsection{Christoffel function weighted least-squares}

Now, we introduce the weighted least-squares procedure for computing the expansion coefficients $\{c_j\}_{j=1}^N$ in the following expansion
\begin{align*}
  f(\xi) \approx \sum_{j=1}^N c_j \Phi_j(\xi),
\end{align*}
with $\{\Phi_j(\xi)\}_{j=1}^N$ being the data-driven orthogonal bases constructed in Section 3, and we denote the associated polynomial space by
\begin{align*}
P_N :=\textmd{span}\Big\{\Phi_j(\xi), \,\, 1\leq j \leq N\Big\}.
\end{align*}
We recall that the polynomial space we considered in this work is of total degree type (\ref{eq:TD}), and the associated maximum polynomial order is denoted by $k.$  The weighted least-squares approach suggest to compute the coefficients via sample evaluations. To this end, suppose we have some sample evaluations $\{f(\mathbf{z}_m)\}$ at some properly chosen samples $\{\mathbf{z}_m\}_{m=1}^M$.
Then, we seek the following weighted discrete least-square approximation $f_N \in P_N$ by requiring
\begin{align}\label{eq:least}
f_N := P^N_m f = \argmin_{p\in P_N} \frac{1}{M} \sum_{m=1}^M \mathbf{w}_m\Big(p(\mathbf{z}_m)-f(\mathbf{z}_m)\Big)^2.
\end{align}
Here $\{\mathbf{w}_m\}_{m=1}^M$ are properly designed weights. An equivalent algebraic formula for the above problem yields:
\begin{align}\label{eq:le}
\mathbf{c}=\argmin_{\mathbf{c}\in \mathbb{R}^{N}} \left\|\mathbf{W^{\frac{1}{2}}A}\mathbf{c}-\mathbf{W^{\frac{1}{2}}f}\right\|^2_2,
\end{align}
where
\begin{align*}
 \mathbf{f}=\big(f(\mathbf{z}_1\big), ...,  f{\mathbf{z}}_m)), \quad \mathbf{A}=\big[\Phi_j(\mathbf{z}_m)\big]\in \mathbb{R}^{M\times N}, \,\,\, j=1,...,N,\,\,\, m=1,...,M,
\end{align*}
and $\mathbf{W}=\textmd{diag}(\mathbf{w}_1, ... ,\mathbf{w}_M)$ is the preconditioning matrix. Notice that in the above approach, the sampling strategy and the pre-conditioner are two key points. Here we shall adopt the strategy in \cite{Narayan_2017CSALS}: Christoffel function weighted least-squares. To this end, we define the associated (scaled) Christoffel-type function of $P_N$ by
\begin{align}\label{eq:Christoffel}
  K(\xi) = \frac{N}{\sum_{j=1}^N \Phi_j^2(\xi)},
\end{align}
The components of the preconditioning matrix $\mathbf{W}$ in our weighted least-squares are evaluations of the (scaled) Christoffel function. i.e.,
\begin{align*}
  \mathbf{w}_m = \frac{N}{\sum_{j=1}^N \Phi_j^2(\mathbf{z}_m)}, \quad m=1,...,M.
\end{align*}
Now, we are ready to  summarize the procedures of our Christoffel function weighted least-squares (more detailed discussions for the sampling strategy and the theoretical motivations will be given later):
\begin{itemize}
\item sampling with respect to the probability density $\widehat{\rho}$ of an equilibrium measure, which depends on the input density $\rho$. When $\xi$ is a random vector with unbounded state space, then $\widehat{\rho}$ also depends on $k$, the maximum polynomial degree of the total degree polynomial space $P_N.$ In this case, we denote the sampling measure by $\widehat{\rho}_k.$

\item evaluate the function $f(\xi)$ (the underlying model) at the selected samples $\{\mathbf{z}_m\}_{m=1}^M.$

\item form $M \times N $ Vandermonde-like matrix $\mathbf{A}$ with entries $\Phi_{n}(\mathbf{z}_m).$

\item form the diagonal preconditioning matrix $\mathbf{W}$ using evaluations of the (scaled) Christoffel function.

\item solve the preconditioned least-squares problem \eqref{eq:le} to approximate the expansion coefficients $\{c_j\}_{j=1}^N$.
\end{itemize}
Notice that in our weighted least-squares approach, the main feature is that the sampling strategy is independent of the data-driven bases, thus this procedure and the associated model simulations can be done in prior. Moreover, we shall show in the following that the sampling strategies are straightforward.

\subsubsection{Sampling measure for bounded domain}

We first consider the bounded case, where we assume (without loss of generality) that the computational domain for $\xi$ is $[-1,1]^d.$  In this case, our sampling measure is always the tensor-product Chebyshev measure, i.e,
$$ \widehat{\rho}(\xi) \sim \frac{1}{\pi^d\prod_{k=1}^d \sqrt{1-\xi_k^2}},$$
regardless of the underling measure (if exists, yet unknown) of the random vector $\xi.$ In other words, the only information we require is that the random variable is located in a bounded domain.

Notice that the equilibrium measure for a bounded domain with \textit{any} admissible input density is the Chebyshev measure, or in other words, the Chebyshev measure is \textit{universal} in the bounded setting. Notice that sampling with Chebyshev measure is straightforward: one can simply generate uniform distributed samples $\{\mathbf{u}_m\}_{m=1}^M$ and then generate $\{\mathbf{z}_m\}_{m=1}^M$ by requiring $$\mathbf{z}_m=\cos(\mathbf{u}_m), \quad m=1,...,M.$$

\subsubsection{Sampling measure for unbounded domain}

We now consider the unbounded case. We remark that very few results are known for the equilibrium measure in unbounded domains. Thus, the results in what follows are our conjectures for which the effectiveness have been well studied numerically in \cite{Narayan_2017CSALS}.

{\textbf{The domain $\mathbb{R}^d$ with Gaussian density}.
We consider the domain $\R^d$ with Gaussian-type input $N(\sigma,\mu)$ (yet the parameters $\sigma,\mu$ can be arbitrary/unknown). As in our setting, we assume that we only have some sample locations, we shall first compute an approximated pair ($\widehat{\mu}, \widehat{\sigma}$) of the input. Then, a simple linear transformation $\widehat{\xi} = (\xi-\widehat{\mu})/\widehat{\sigma}$ can be used to make sure that the input has distribution  $N(1,0)$ (approximated).  A conjecture result \cite{Narayan_2017CSALS, Jakeman_2016generalizedsample} for the equilibrium measure associated with $N(1,0)$ is given by
\begin{align*}
  \widehat{\rho}(\xi) = C \left( 2 - \left\|\xi\right\|^2 \right)^{d/2},
\end{align*}
with $C$ a normalization constant. Furthermore, we shall expand the associated samples (generated by the above measure) by the square root of the maximum polynomial degree $k$. The following is a concrete way to sample from this expanded density:
\begin{enumerate}
  \item Compute $k$, the maximum polynomial degree of $P_N$.
  \item Generate a vector $\mathbf{y}=(y_1,\ldots,y_d)$ of $d$ independent normally distributed random variables.
  \item Draw a scalar sample $\nu$ from the Beta distribution on $[0,1]$, with distribution parameters $\alpha=d/2$ and $\beta=d/2+1$.
  \item Finally, we set
    \begin{align*}
        \mathbf{z}= \frac{\mathbf{y}}{\|\mathbf{y}\|_2}  (2k\nu)^{\frac{1}{2}}.
    \end{align*}
\end{enumerate}
The above procedure generates samples on the Euclidean ball of radius $\sqrt{2 k}$ in $\R^d$. We emphasize that our methodology samples from a density that is only a conjecture for the correct equilibrium measure. We also remark that we have introduced a density error to this approach, as the mean and variance are computed approximately. How to quantify and control such errors will be our future projects.

\textbf{The domain $\mathbb{R}^d_+$ with exponential density.} Let $\xi$ take values on $\mathbb{R}^d_+$ with associated exponential-type probability density (again the associated parameters can be arbitrary). Again, we shall compute an approximated mean value so that we can work with the standard exponential-type probability density. In this case we sample from the following density function
\begin{align*}
  \widehat{\rho}(\xi) = C \sqrt{\frac{\left(4 - \sum_{i=1}^d \xi_i\right)^d}{\prod_{i=1}^d \xi_i}}
\end{align*}
As we conjectured in \cite{Narayan_2017CSALS, Jakeman_2016generalizedsample}, this is the equilibrium measure associated to this choice of $\rho$. We shall also expand the samples by the maximum polynomial degree $k.$ The following is a concrete way to sample from this expanded density:
\begin{enumerate}
  \item Compute $k$, the maximum polynomial degree of the polynomial space.
  \item Generate a $(d+1)$-dimensional Dirichlet random vector $\mathbf{y}$ with parameters
  $\left(\frac{1}{2}, \frac{1}{2}, \ldots, \frac{1}{2}, \frac{d}{2}+1\right)$.
  \item Truncate the last ($(d+1)$'th) entry of $\mathbf{y}.$
  \item Set $ \mathbf{z}= 4 k \mathbf{y} $.
\end{enumerate}
\begin{Remark}
In the above, we have only discussed two most commonly used densities in unbounded domains, i.e., the Gaussian density and the exponential density. For more general unbounded densities, less is known for the equilibrium measure (even in the conjecture sense). A possible way to handle such situations is to truncate the domain into a finite one, and then perform the Chebyshev sampling in the finite domain. However, this is non-trivial due to the truncated error and we left such cases for future studies,
\end{Remark}

\subsection{Theoretical motivations}

In this section, we shall provide with some motivations for our Christoffel weighted least-squares. We shall only show the motivation in the bounded domain setting, and one can refer to \cite{Narayan_2017CSALS} for the motivation of unbounded domain cases. To begin, we first present the following fundamental result for the least-squares stability \cite{Cohen_LS}:
\begin{theorem}
For a $d$-dimensional function $f(\xi)$, consider its approximation in a finite orthogonal bases space $P_N=\textmd{span}\{\Phi_j(\xi), \,\, 1\leq j \leq N\}$ with the associated orthogonal density $\rho(\xi).$ Suppose the samples $\{\mathbf{z}_m\}_{m=1}^M$ are generated with respect to $\rho(\xi).$ Consider the following least-squares approach
\begin{align}
\mathbf{c}=\argmin_{\mathbf{c}\in \mathbb{R}^{N}} \left\|\mathbf{A}\mathbf{c}-\mathbf{f}\right\|^2_2,
\end{align}
Then, the above algorithm is stable in the following sense
\begin{align*}
\mathbf{Pr}\left\{ \| \mathbf{A-I}\| \geq \frac{1}{2} \right\} \leq 2 M^{-r}
\end{align*}
provided that
\begin{align*}
\kappa(N):= \max_{\xi} \sum_{j=1}^N \Phi^2_j(\xi) \leq \delta \frac{M}{\log M} \quad \textmd{with} \quad  \delta= \frac{1-\log 2}{2-2r}.
\end{align*}
Here $\mathbf{I}$ is the identity matrix.
\end{theorem}

The above theorem states that to make the algorithm stable, it is essential to control the quantity $\kappa(N)$ as one requires approximately $M \gtrsim \kappa(N)$ (up to a logarithmic factor). However for many cases, the quantities $\kappa(N)$ behaves super-linear in $N$ leading to too much demanding conditions on the sampling size
$M$ to guarantee stability. For example, the most commonly used Legendre polynomials gives $\kappa(N)\thicksim N^2$ meaning that one requires $M \geq C N^2,$ which is not satisfactory.

The above observations motivate us to use a weighted version of least-squares. In our approach, by introducing the pre-conditioner $\mathbf{W},$ we are in fact working with a scaled bases set (see (\ref{eq:Christoffel}) for the definition of $K(\xi)$)
\begin{align}\label{eq:Q-def}
  \widehat{P}_N = \mathrm{span} \left\{ \widehat{\Phi}_j= \frac{\Phi_j}{\sqrt{K(\xi)}} \;\; \big| \;\; 1\leq j \leq N \right\}.
\end{align}
It is easy to show that for the new bases $\widehat{\Phi}_j$ it holds
\begin{align}
\widehat{\kappa}(N):= \max_{\xi} \sum_{j=1}^N \widehat{\Phi}^2_j(\xi)  \equiv N.
\end{align}
This means that we have the optimal control of the associated quantity $\widehat{\kappa}(N)$.

However, to show the optimal stability by Theorem 4.1 (which use samples according to the orthogonal measure), we have to sampling with a transformed measure
\begin{align}
\widetilde{\rho}(\xi) \sim  K(\xi) \rho(\xi) = \frac{N \rho(\xi)}{\sum_{j=1}^N \Phi^2_j(\xi)},
\end{align}
as our new bases are orthogonal according to $\widetilde{\rho}(\xi)$. Notice that $\widetilde{\rho}(\xi)$ depends on the polynomial space, and furthermore, sampling with $\widetilde{\rho}(\xi)$  seems to be non-trivial. Nevertheless, we learn from potential theory in the bounded setting that \cite{Narayan_2017CSALS}
\begin{align}
\widetilde{\rho}(\xi) \rightarrow \widehat{\rho}(\xi),  \quad \textmd{when} \quad  N\rightarrow \infty.
\end{align}
The above result motivated us to sampling with $\widehat{\rho}(\xi)$ -- the equilibrium measure. In this way, we can get a stable approach in the asymptotical sense ($N\rightarrow \infty$).  And furthermore, our sample strategy now is independent of the polynomial space, and this is advantage for adaptive computations where the polynomial spaces are constructed adaptively.
In the bounded setting, for any admissible input density, the equilibrium measure is just the Chebyshev density, and this is the exact motivation for us to introduce the Christoffel weighted least-squares.

\section{Numerical experiments}
In this section, we present several numerical examples to show the effectiveness of our Christoffel weighted least-squares for data-driven polynomial approximations. We are interested primarily in investigating how the sampling rates between $M$ and $N$ affect stability and accuracy. Due to the probabilistic nature of the random sampling method, all reported results are averaged over 100 independent tests to reduce the statistical oscillations. In all our figures and numerical tests, we shall show the performance with a linear and a $\log$-linear dependence between $M$ and $N,$ namely, $M = CN$ and $M= CN \log{N}.$  The following stochastic input distributions will be considered:
\begin{itemize}
\item Discrete Binomial distribution: Bino(n,p) in $[-1,1]$:
\begin{align*}
f(k;n,p)=\mathcal{P}(\xi=\frac{2k}{n}-1)=\frac{n!}{k!(n-k)!}p^k(1-p)^{n-k}, \ k=0,1,\ldots,n;
\end{align*}

\item Discrete Poisson distribution: Pois($\lambda$) in $[-1,1]$:
\begin{align*}
f(k|\lambda)=\frac{\lambda^k}{k!}\exp(-\lambda);
\end{align*}

\item Uniform distribution: $U[a,b]:$
\begin{align*}
f(x)=\begin{cases}
\frac{1}{b-a}, & x\in[a,b]\\
0, & \text{otherwise}.
\end{cases}
\end{align*}

\item Exponential distribution $\textmd{Exp}(\mu)$ in $(0, \infty)$ with parameters $\mu$:
\begin{equation*}
f(x|\mu)=\frac{1}{\mu}\exp\big(-\frac{x}{\mu}\big).
\end{equation*}

\item Normal distribution $N(\mu,\sigma)$ in $(-\infty, \infty)$ with parameters $\mu, \sigma$:
\begin{equation*}
f(x|\mu,\sigma)=\frac{1}{\sqrt{2\pi}\sigma}\exp\big(-\frac{(x-\mu)^2}{2\sigma^2}\big).
\end{equation*}
\end{itemize}

\begin{table}[!ht]\centering
\caption{Test examples for the two-dimensional case.}
\label{tab:2dRV}
\begin{tabular}{c|l}
\text{Type} &  \text{Parametric Distributions} \\ \hline
1 & $\xi_1\sim \textmd{Bino}(20,1/2), \,\,\, \xi_2\sim U[-0.6,0.6]$ \\ \\
2 &  $\xi_1\sim U[-0.8,0.8], \,\,\,\xi_2\sim U[-1,1]$\\ \\
3 &  $\xi_1\sim \textmd{Bino}(20,1/2),\,\,\, \xi_2\sim \textmd{Pois}(10)$\\ \\
4 &  $\xi_1\sim  U[-0.6,0.6], \,\,\,\, \xi_2\sim N(0.1,1.2)$\\ \hline
\end{tabular}
\end{table}

\subsection{Stability tests}
We first test the condition number of the design matrix
$$\textmd{Cond}(\mathbf{\widehat{A}}) = \frac{\lambda_{max}(\mathbf{\widehat{A}})}{\lambda_{min}(\mathbf{\widehat{A}})} \qquad \textmd{with} \qquad  \mathbf{\widehat{A}} =\mathbf{W^{\frac{1}{2}}A}.$$ The main focus is how this quantity is affected by the the sampling rate $M/N$. Notice that this quantity measures the sensitivity of the solution of a system of linear equations to errors in the data, that is, it directly reflects the stability of the method. In all examples that follow we perform 100 trials of each procedure and report the mean condition number along with $20\%$ and $80\%$ quantiles.

\begin{figure}[htbp]
\begin{center}
    \includegraphics[width=6cm]{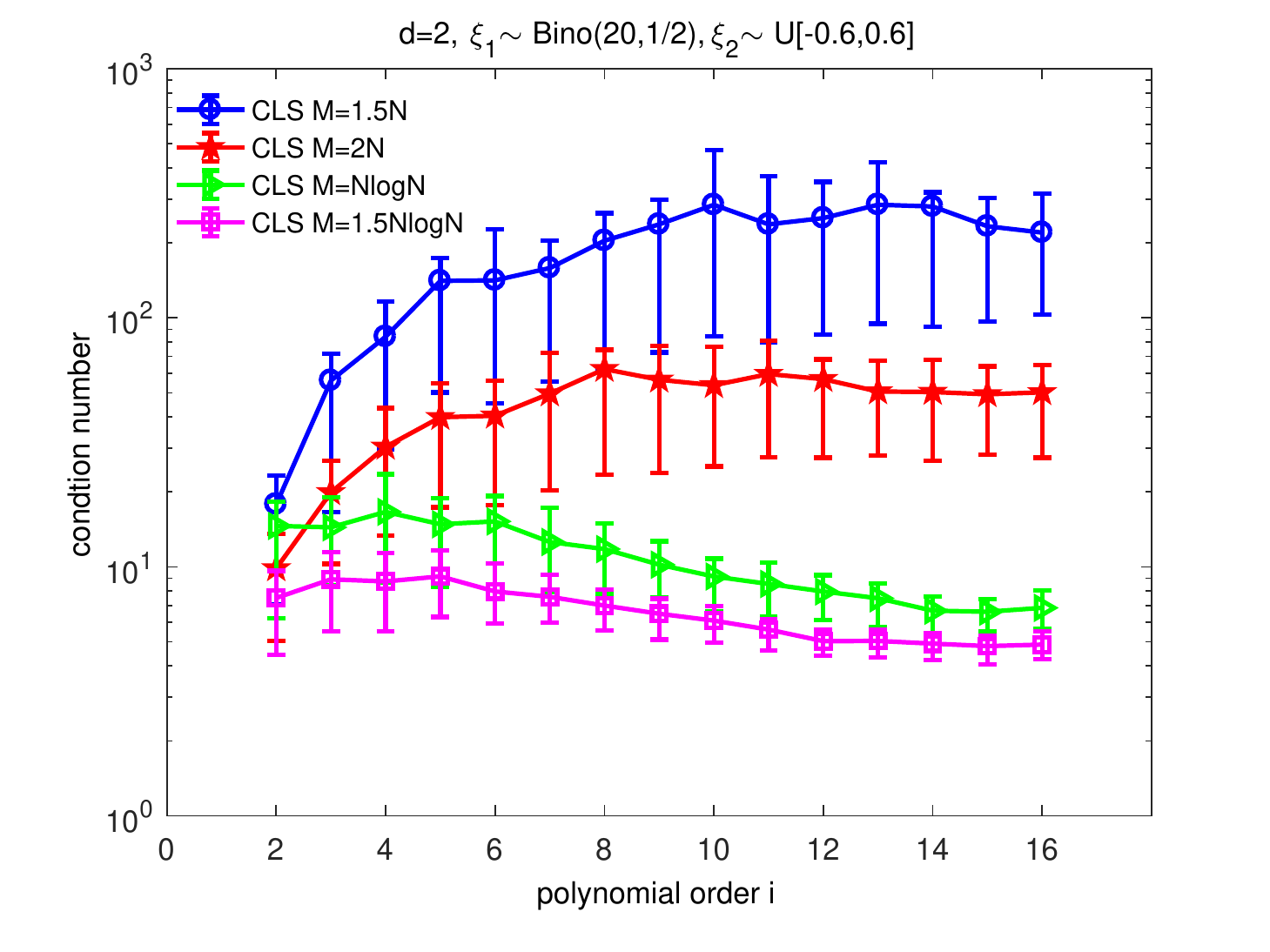}
    \includegraphics[width=6cm]{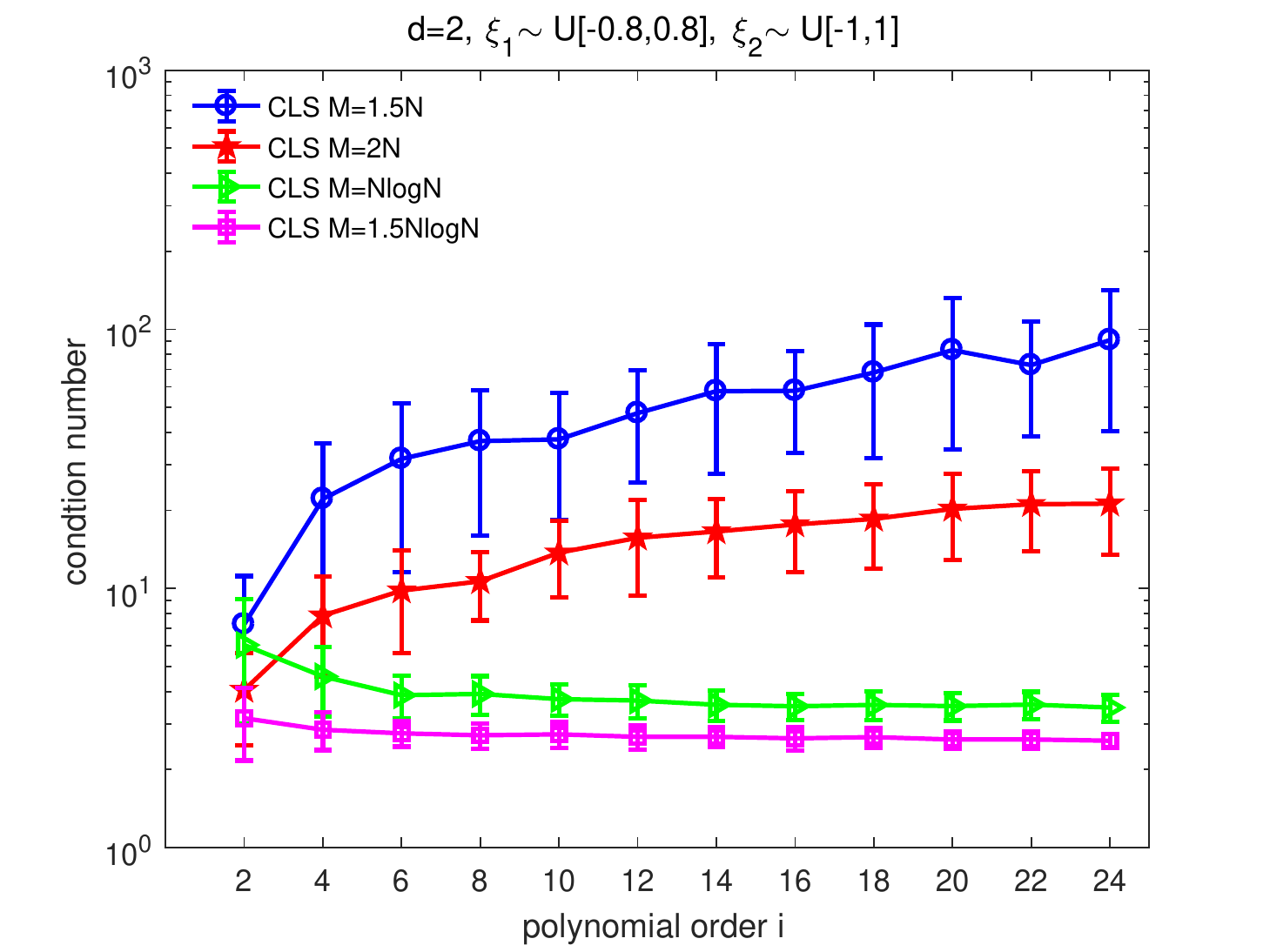}
    \includegraphics[width=6cm]{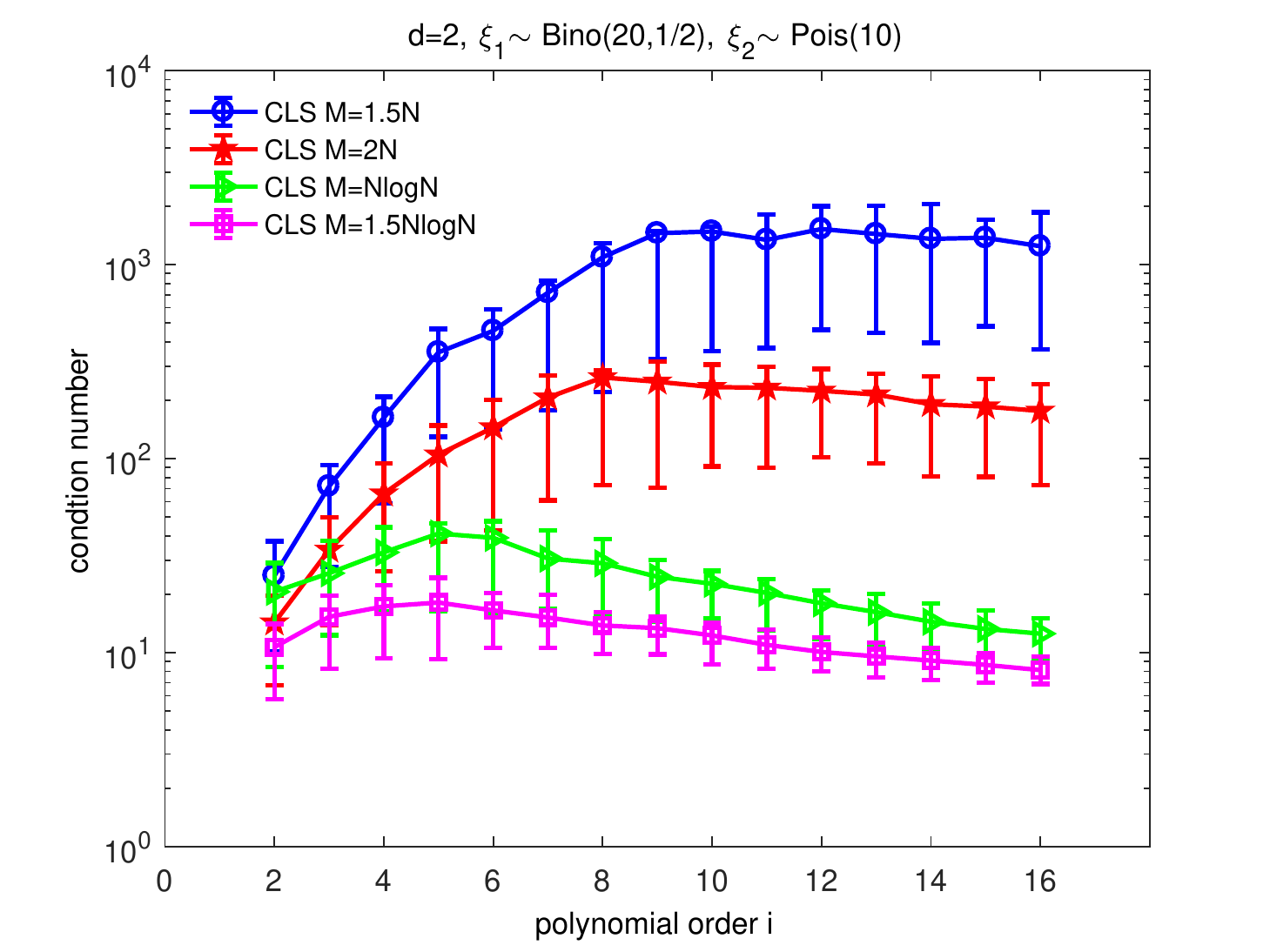}
     \includegraphics[width=6cm]{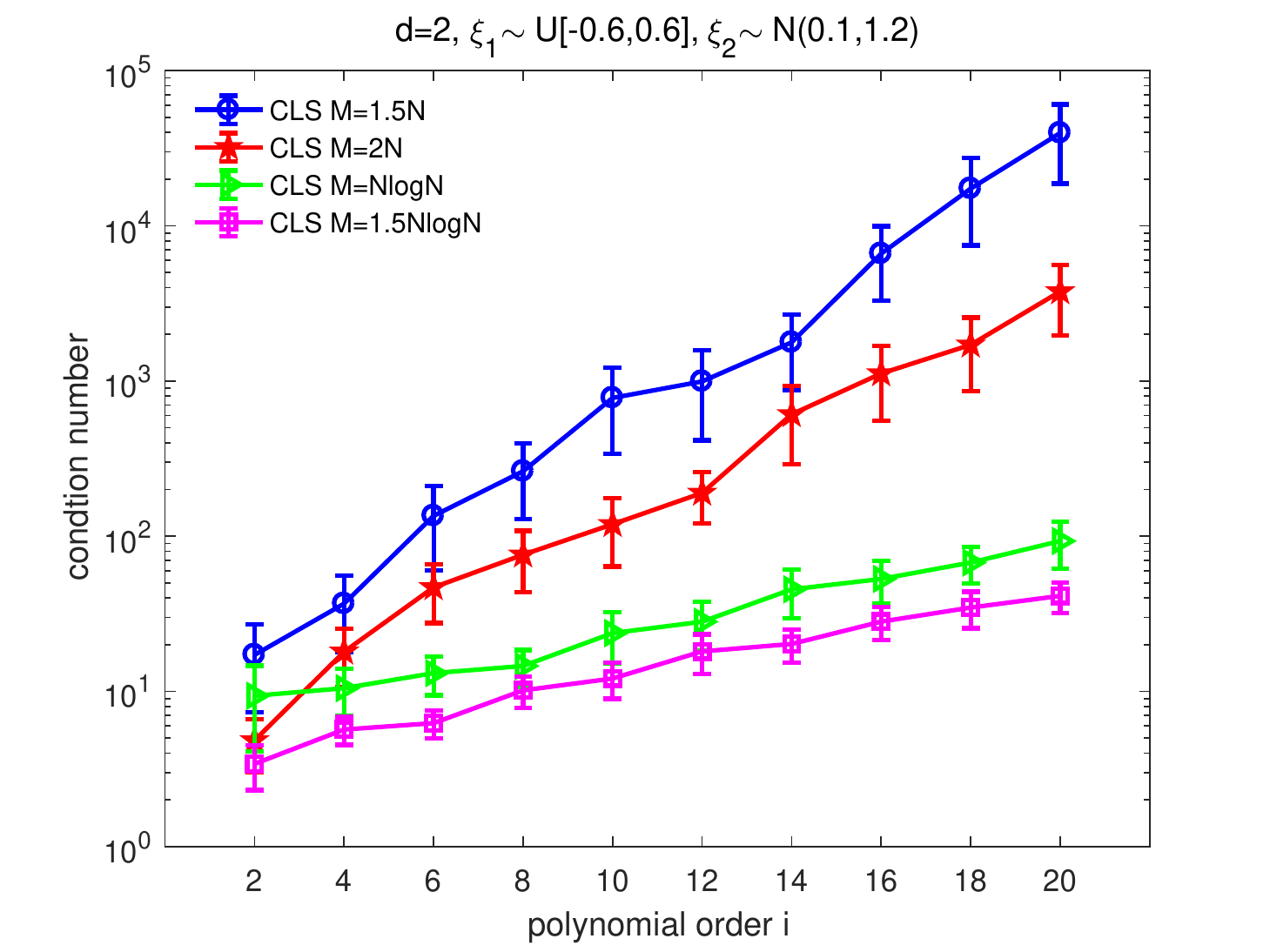}
\end{center}
  \caption{Condition number with respect to the polynomial degree in the 2-dimensional case (Table \ref{tab:2dRV}) with different sampling rates.}\label{fig:2dcond}
\end{figure}
We first consider the two dimensional tests. Four different test cases are given in Table \ref{tab:2dRV}, where the uniform distribution with different parameters for each dimension and mixture distributions (including binomial, poisson distribution and normal) are taken into account. Notice that the fourth test case includes both bounded and unbounded distributions, and thus in our test, we shall sampling with different equilibrium measures in each dimension. Here we use the associated moments directly (so that the numerical error for computing the moments with samples can be neglected) to construct the data-driven bases (by the moment match method). Then, we sampling with the equilibrium measure and construct the associated design matrix. We have presented the condition numbers of the design matrix for the four test cases in Fig. \ref{fig:2dcond}. Different sampling rates are reported, i.e, $M=1.5N,$ $M=2N,$  $M=N\log N,$ and $M=1.5N\log N.$ We notice that the $\log$-linear sampling rate produces more stable results -- the condition number is bounded above for the first three test cases. However, for the fourth test case, we still observed a slightly growing trend, and this is due to the involved unbounded random variable.

We next consider a synthetic example for an empirical data distribution. The simulation data set is generated as the superposition of uniform, normal and log-normal distributions (with sample size $M=10000$), see Fig. \ref{fig:2dcond_hist} (Left). Here we construct the data-driven polynomial bases based on the moments that are computed by those samples. The corresponding condition number for this test case is shown in Fig. \ref{fig:2dcond_hist} (Right). Again, we observe that the $\log$-linear sampling rate provides more stable result.
\begin{figure}[htbp]
\begin{center}
\includegraphics[width=6cm]{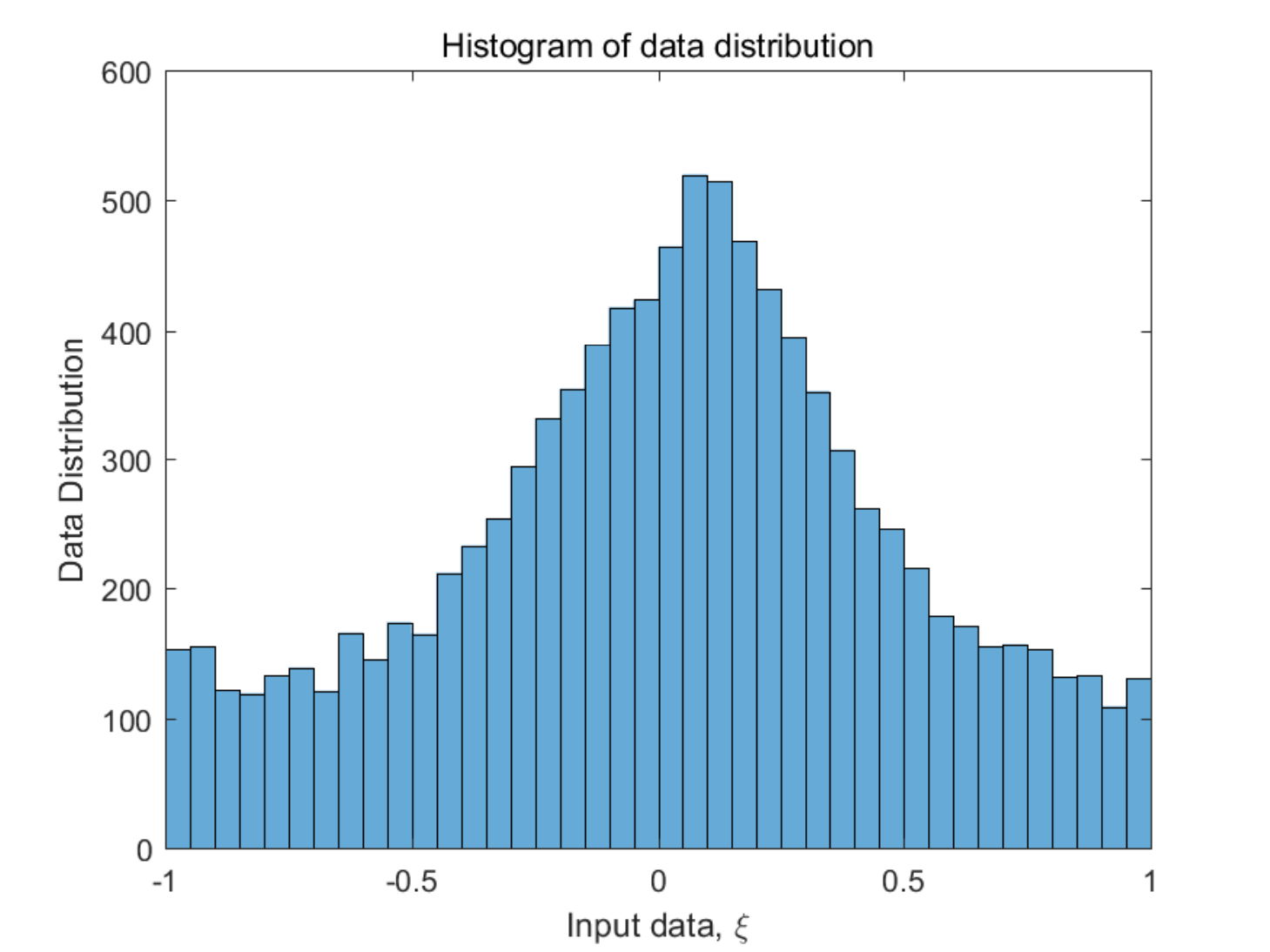}
     \includegraphics[width=6cm]{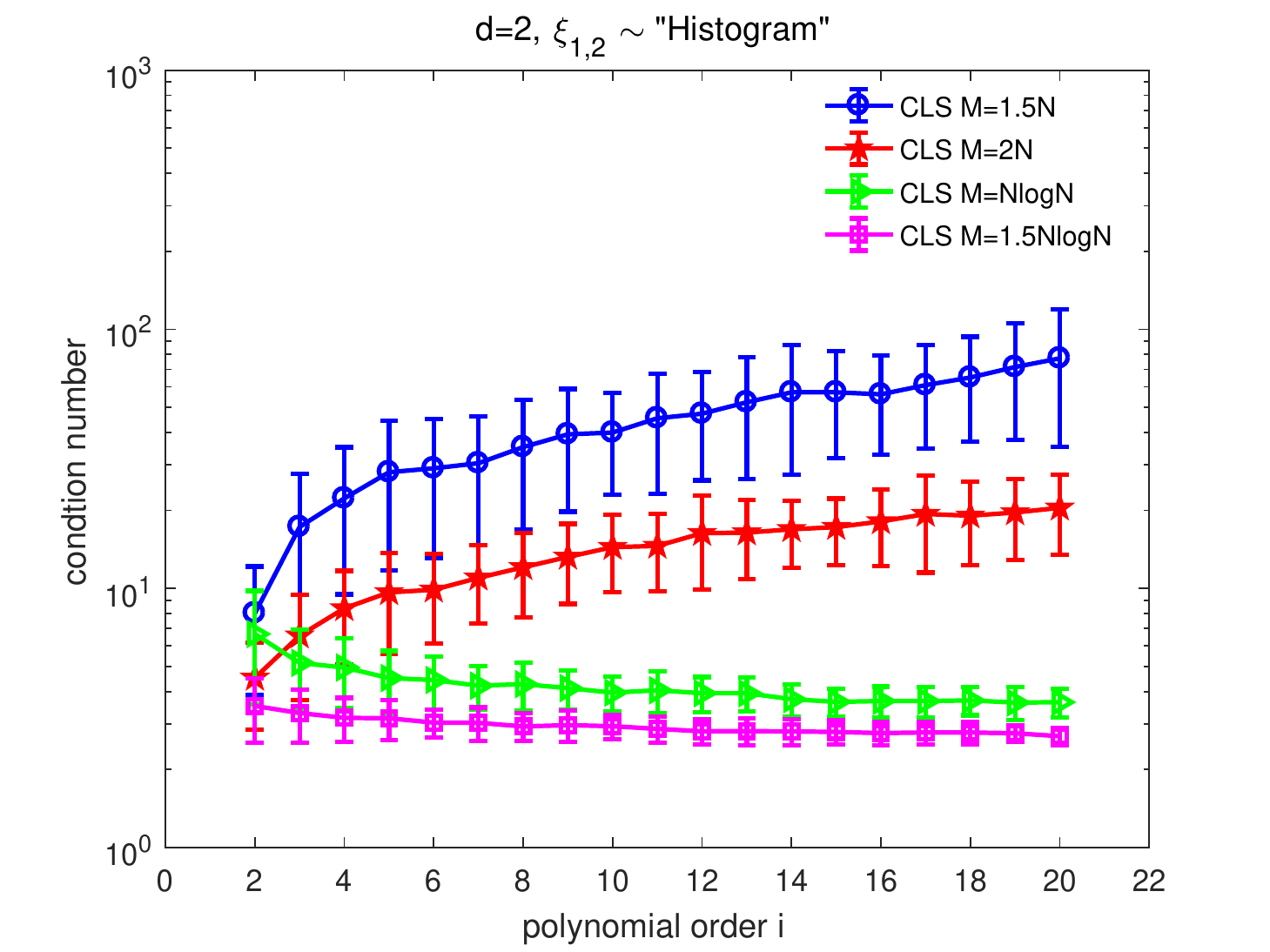}
\end{center}
  \caption{Left: Data distribution showed by histogram. Right: Condition number with respect to the polynomial degree in the 2-dimensional polynomial spaces. }\label{fig:2dcond_hist}
\end{figure}

Finally, we further test the stability for the five dimensional case. The input random parameters used are listed in Table \ref{tab:5dRV}, and the corresponding condition numbers of the different test cases are reported in Fig. \ref{fig:5dcond}. For all test cases, the design matrix admits more stable property with the $\log$-linear sampling rate. However, for cases that involve unbounded parameters, we can still observe a slightly growing trend.

\begin{table}[!ht]\centering
\caption{Test examples for the five dimensional case.}
\label{tab:5dRV}
\begin{tabular}{c|l}
\text{Type} &  \text{Parametric Distributions} \\ \hline
1 & $\xi_i\sim U[a_i,b_i], \,\,\, a=[-0.1,-0.5,-0.8,-1,-1.2], \,\,\, b=-a.$ \\ \\
2 & $\xi_i\sim N(\mu_i,\sigma_i), \,\,\, \mu=[0,01,-0.1,0.2,-0.2], \,\,\, \sigma=[1,1.1,1.2,1,0.9].$\\ \\
3 & $\xi_{1,2}\sim U[-0.6,0.6],  \,\,\, \xi_{3,4}\sim \textmd{Bino}(20,1/2), \,\,\, \xi_5\sim \textmd{Pois}(10).$\\ \\
4 & $\xi_1, \xi_2\sim U[-1,1], \,\,\, \xi_3\sim N(0,1), \,\,\, \xi_4\sim N(0.1,1.5), \,\,\, \xi_5\sim N(0.2,2).$\\ \hline
\end{tabular}
\end{table}

\begin{figure}[htbp]
\begin{center}
     \includegraphics[width=6cm]{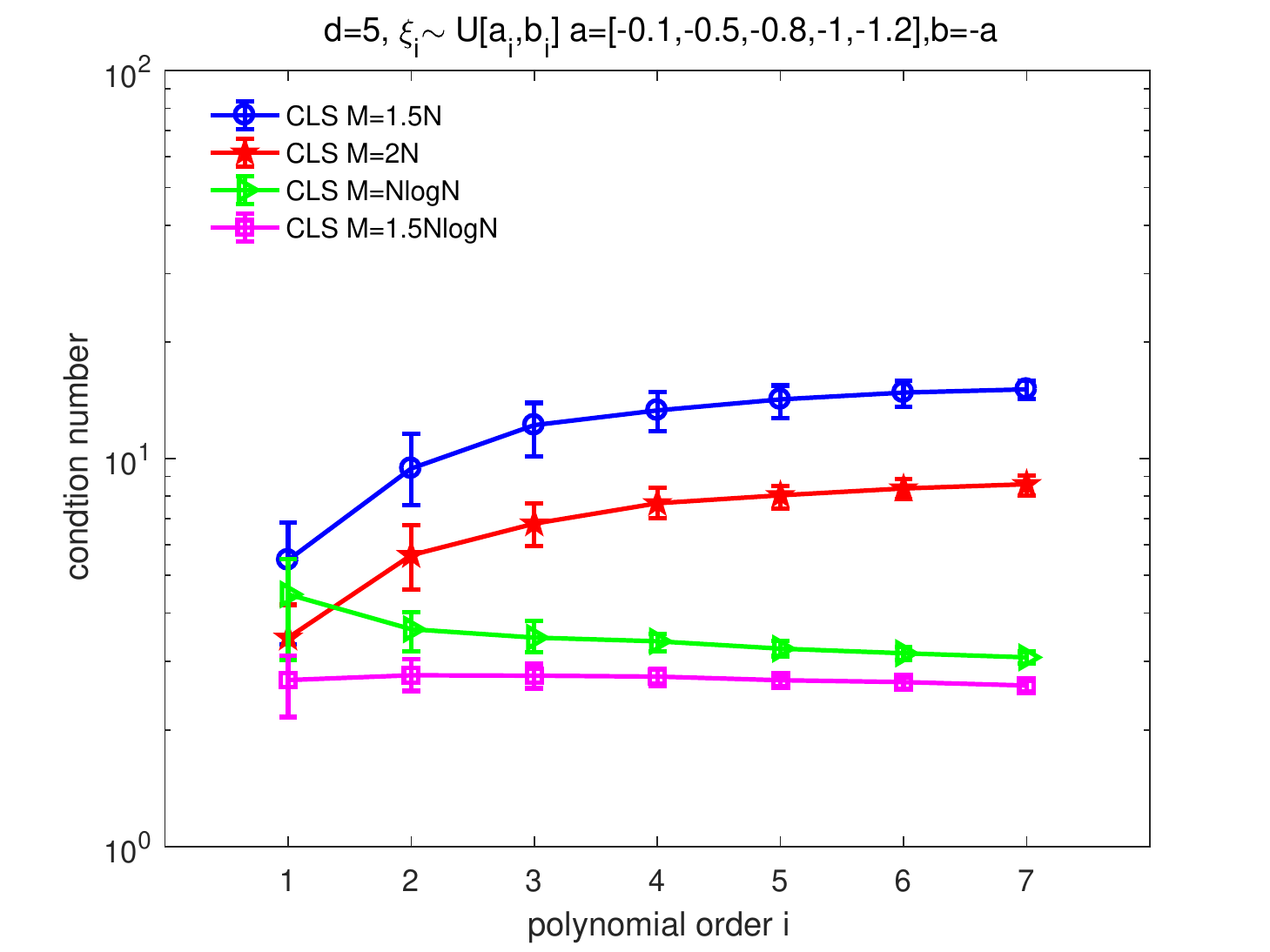}
     \includegraphics[width=6cm]{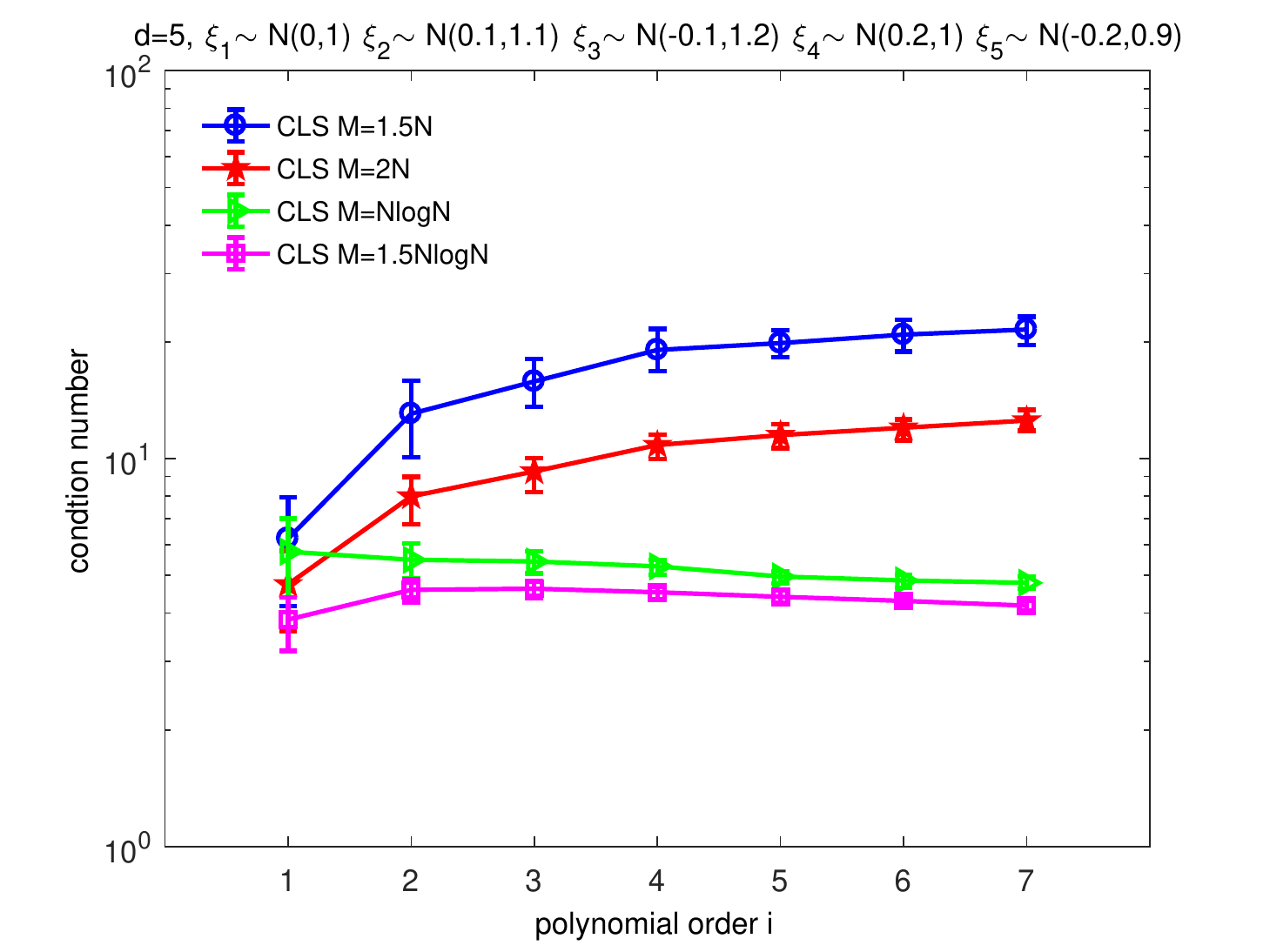}
     \includegraphics[width=6cm]{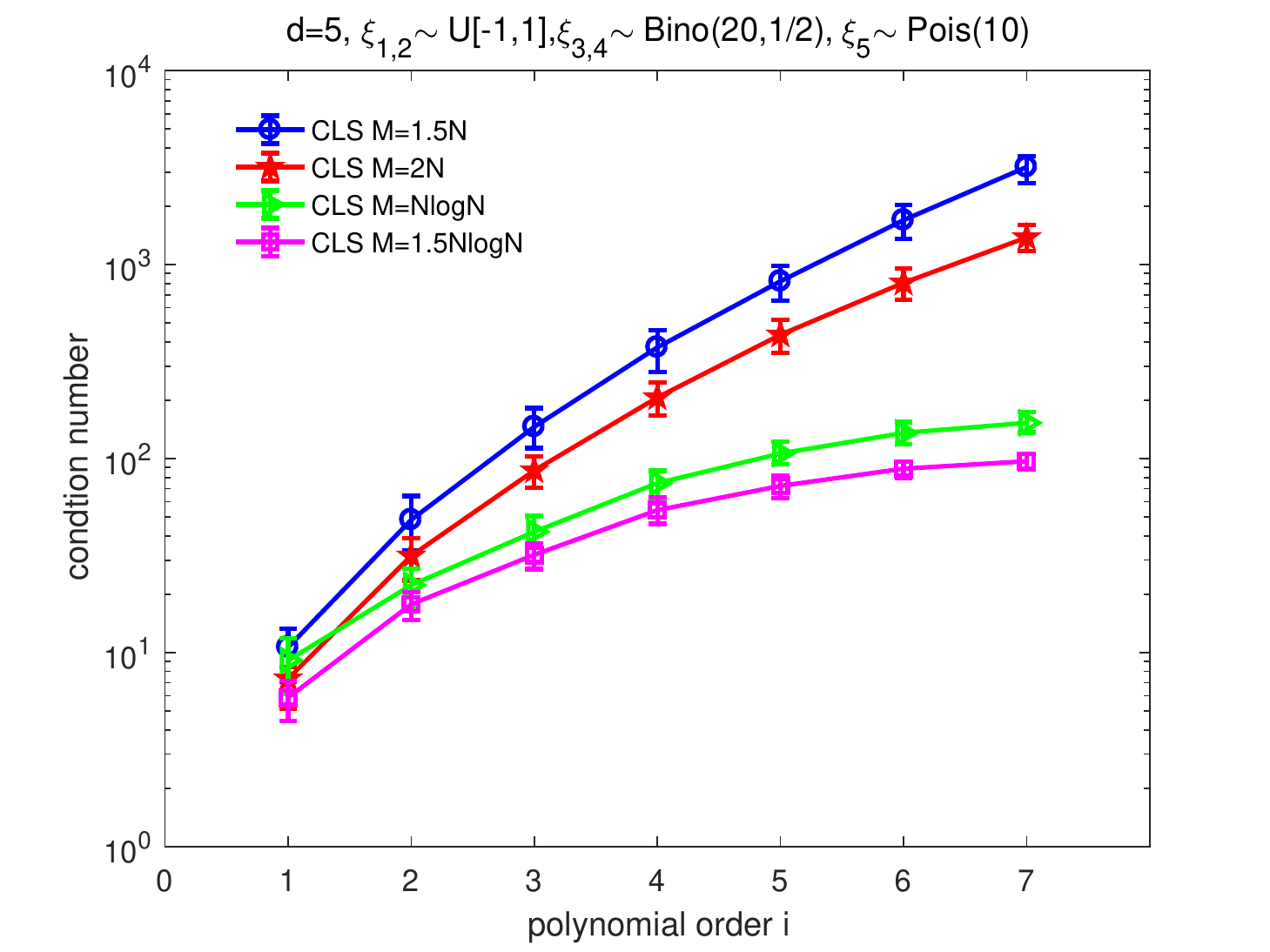}
     \includegraphics[width=6cm]{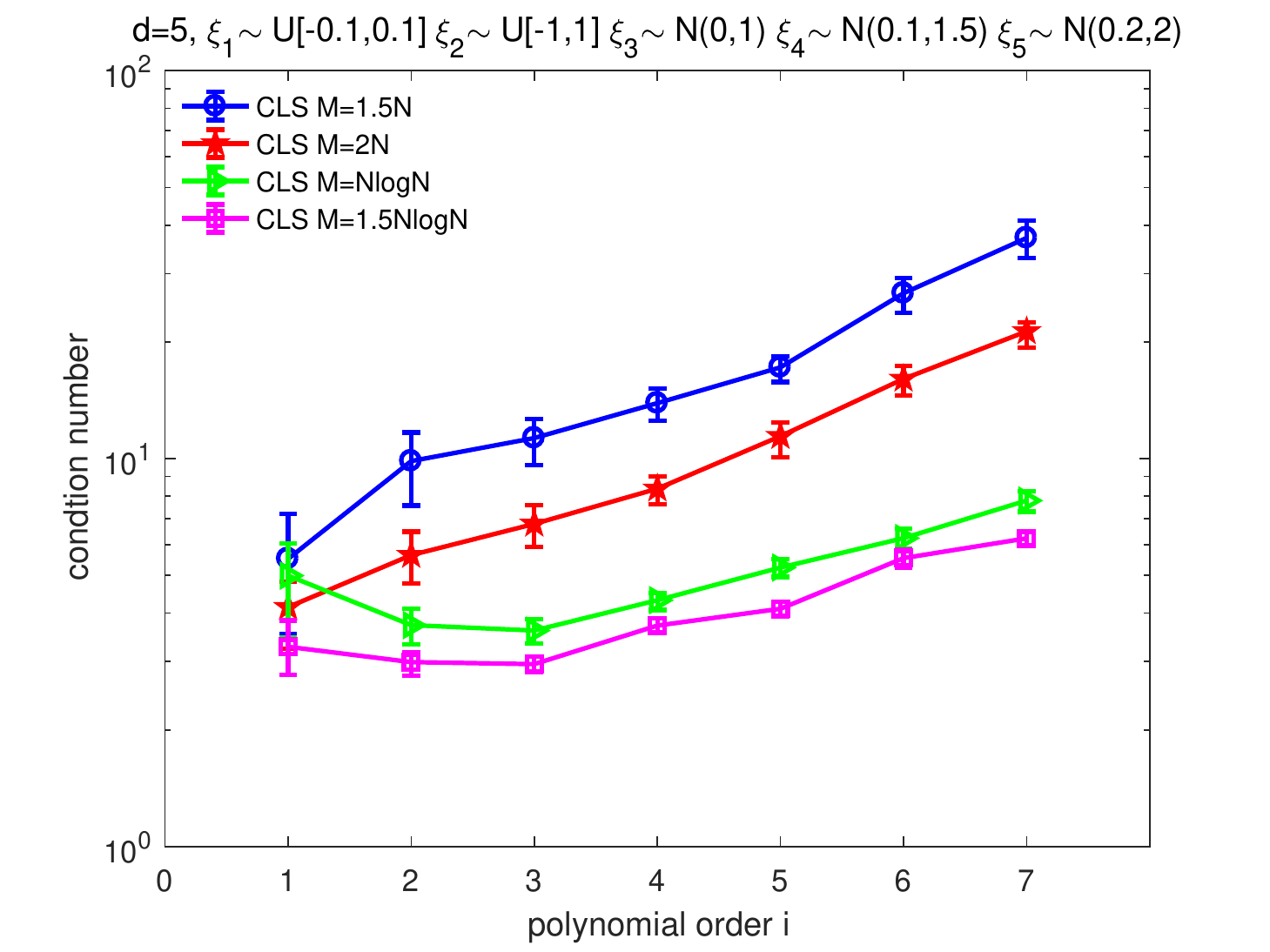}
\end{center}
  \caption{Condition number with respect to the polynomial degree for the five-dimensional tests in Table \ref{tab:5dRV}.}\label{fig:5dcond}
\end{figure}

\subsection{Accuracy tests}
We now test the approximation accuracy of the data-driven bases with Christoffel least-squares post-processing. We shall use the discrete $\ell_2$-error to measure the performance of the approximation, namely, for a given function $f(\xi)$ and a given set of random samples $\{\textbf{z}_l\}_{l=1}^{L}$ in the state space, we evaluate the numerical error via
\begin{align*}
\epsilon =\left(\frac{1}{L}\sum_{l=1}^{L}|f_N(\textbf{z}_l)-f(\textbf{z}_l)|^2\right)^{1/2},
\end{align*}
where $f_N$ is the lease-square solution using the data-driven bases.

\subsubsection{Function approximations}
We first consider the following different test functions:
\begin{align*}
f_1(\xi)=\exp\left(\sum_{k=1}^d\xi_k\right), \quad f_2(\xi)=\sum_{k=1}^d 0.3+\sin\left(\frac{16}{15}(\xi_k-0.7)\right)+\sin^2\left(\frac{16}{15}(\xi_k-0.7)\right)\\
f_3(\xi)=\exp\left(-\sum_{k=1}^d c_k^2(\xi_k-0.01)^2\right), \quad  c_k=\exp\left(-6k/d\right), \quad
f_4(\xi)=\sin\left(\sum_{k=1}^d\xi_k\right).
\end{align*}
The distribution information for the above parameters coincides with Table \ref{tab:2dRV}. The convergence rates of our approach for the two-dimensional case are presented in Fig. \ref{fig:unbounded_f54err_d2}. It is clear shown that the Christoffel least-squares provide very stable and accurate approximation results. In Fig. \ref{fig:err_d5}, we have also tested the five dimensional cases with parameters defined in Table \ref{tab:5dRV} (type 1 and 4), for the test functions $f_1(\xi)$ and $f_3(\xi),$ respectively. Again, our approach admits very stable approximation results.

Finally, we consider tests with histograms data for both the two and five dimensional cases. We consider two sets of data generated as superposition of uniform, normal and log-normal distributions. Results of these approximations are given in Figs. \ref{fig:rawdata_ferr_d2} and \ref{fig:rawdata_ferr_d5}.

\begin{figure}[htbp]
\begin{center}
\includegraphics[width=6cm]{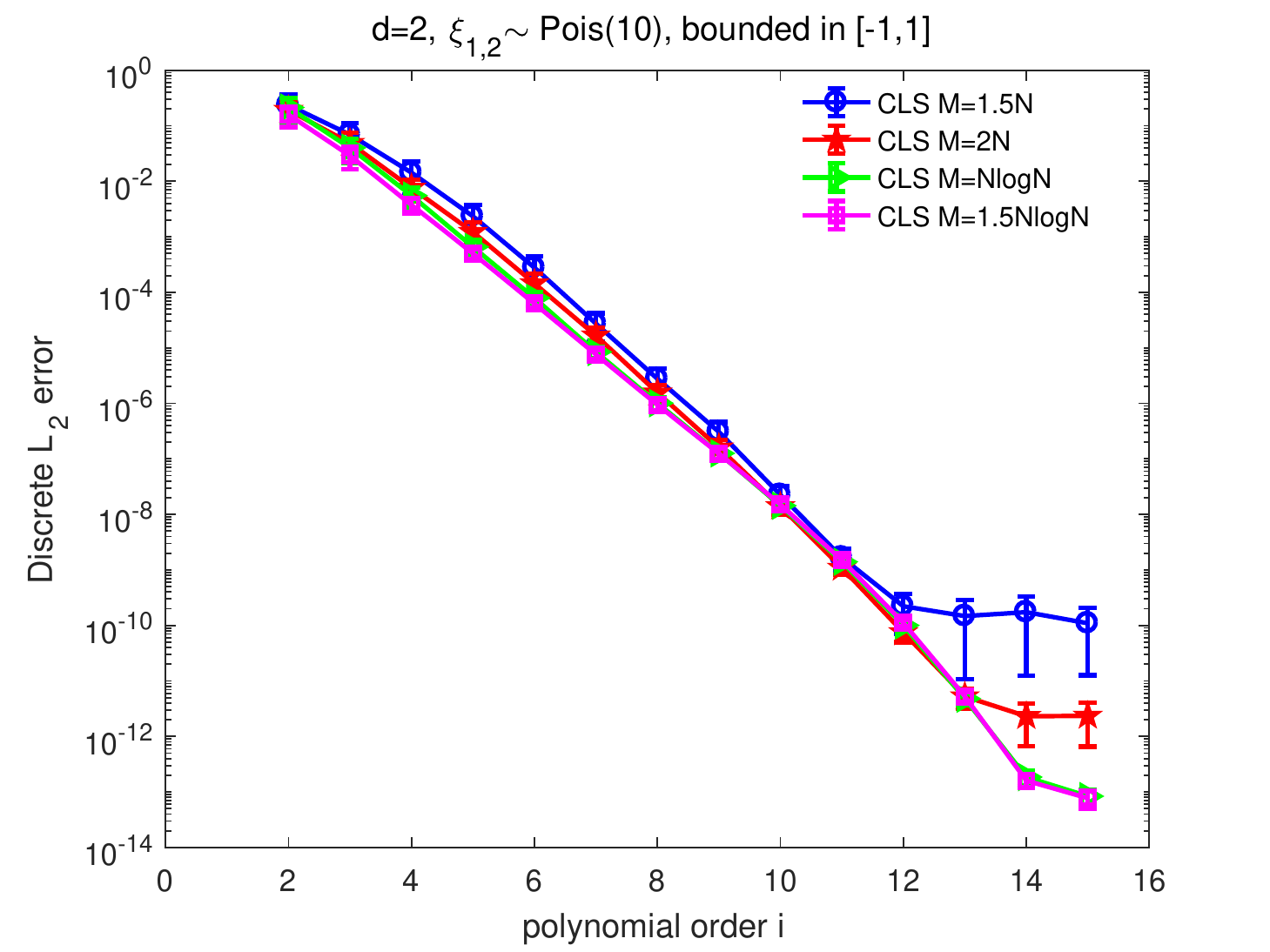}
 \includegraphics[width=6cm]{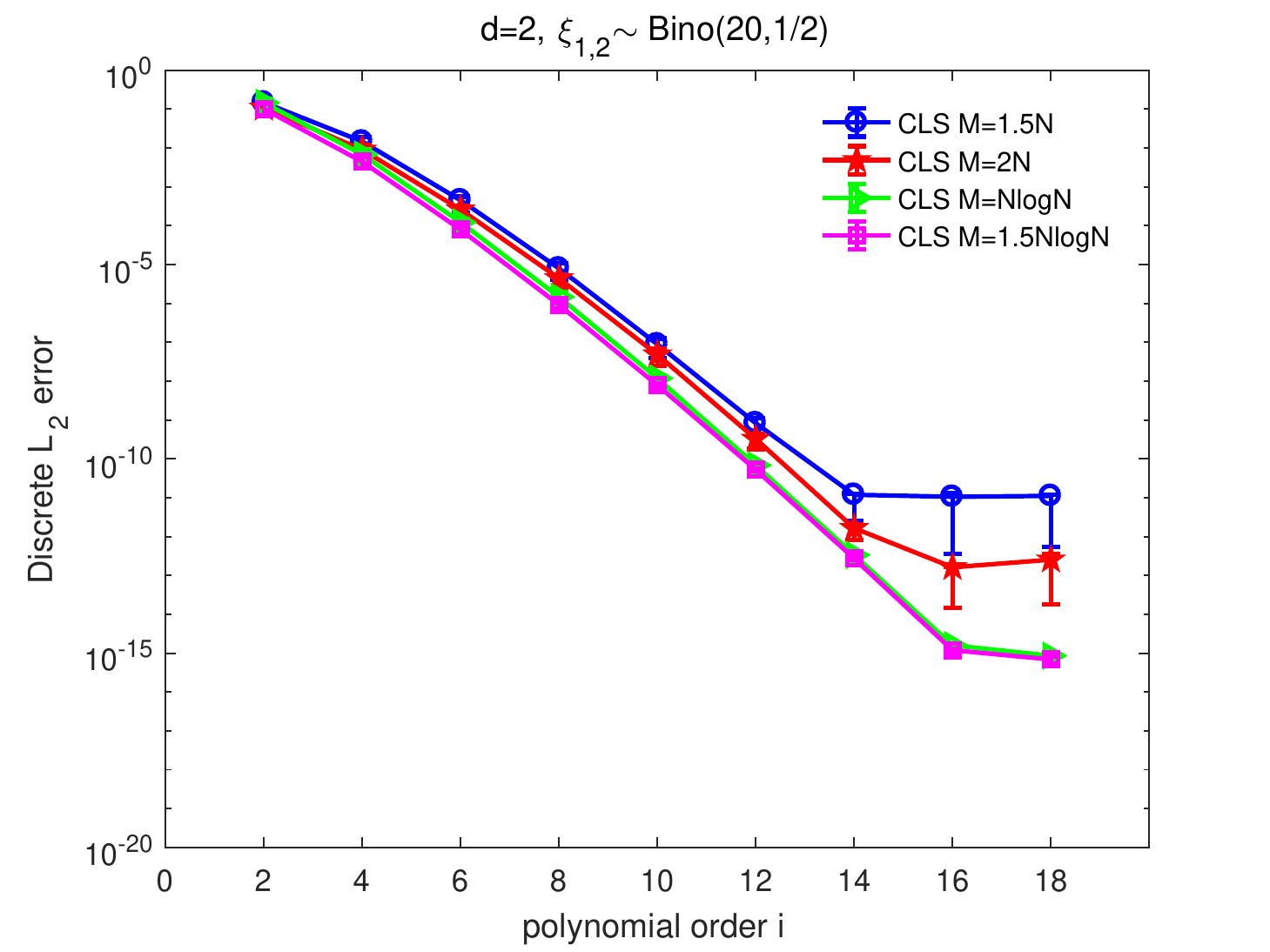}
\\
\includegraphics[width=6cm]{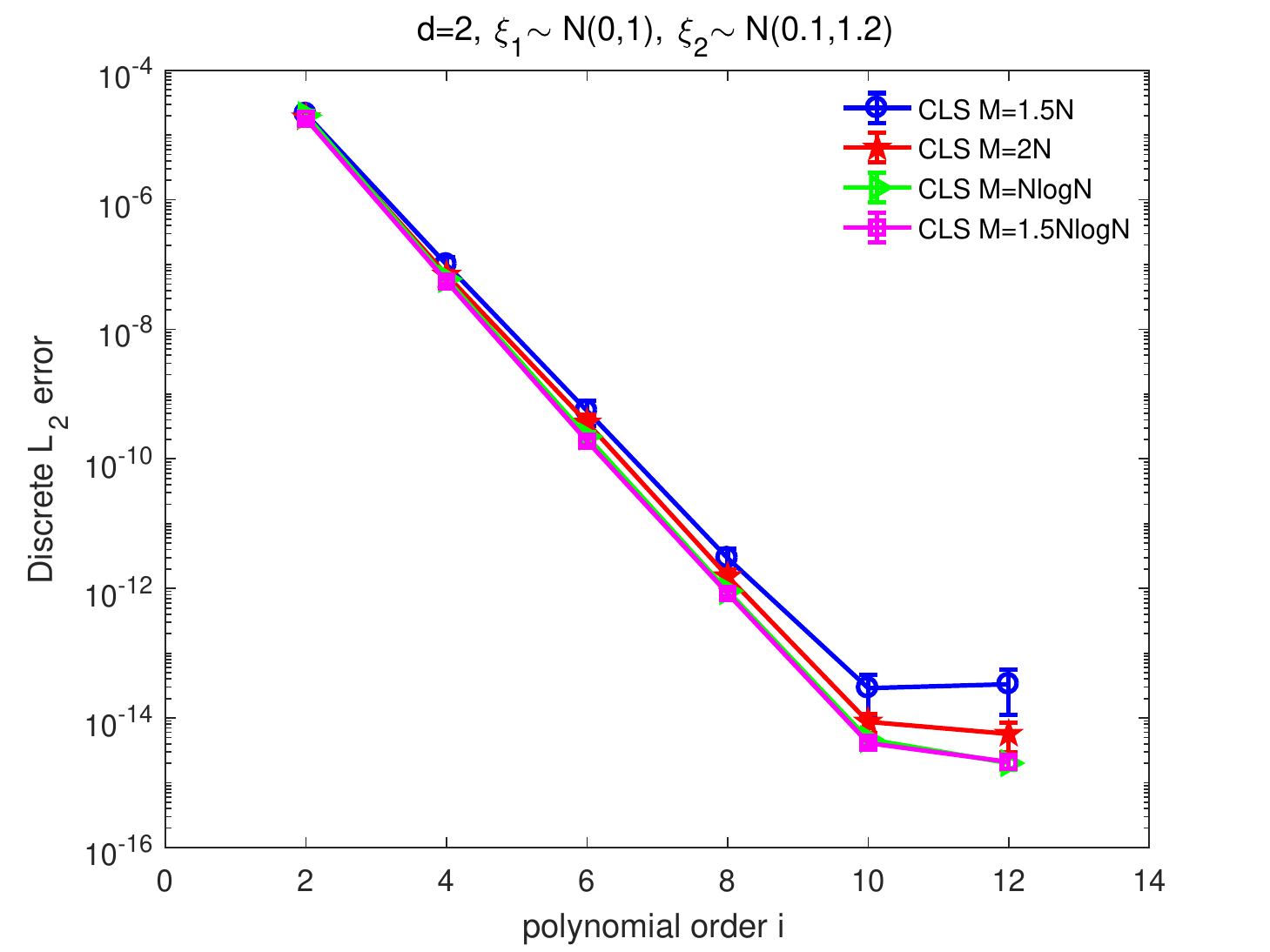}
\includegraphics[width=6cm]{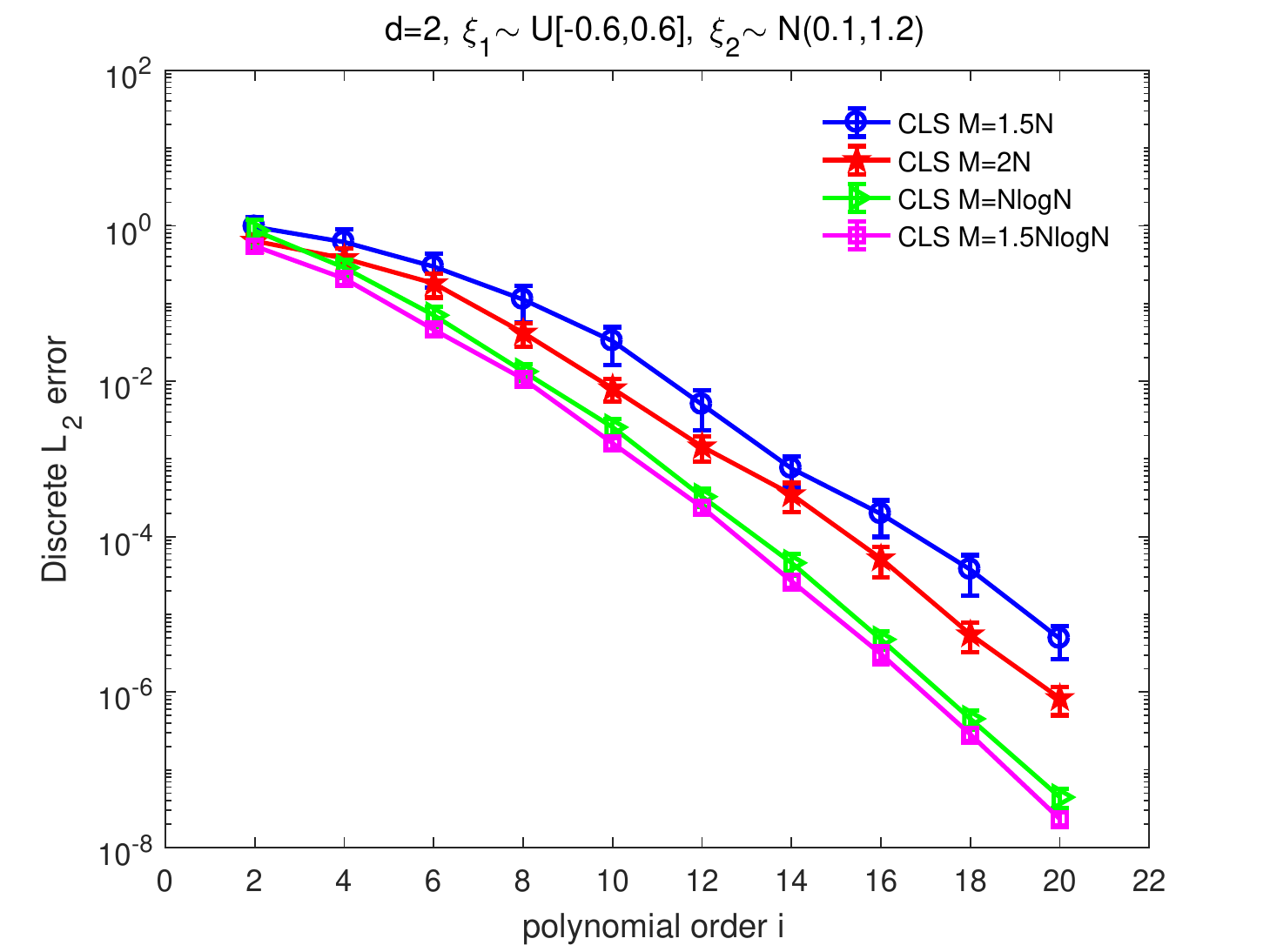}
\end{center}
\caption{Approximation error against polynomial degree for the two dimensional case.}\label{fig:unbounded_f54err_d2}
\end{figure}

\begin{figure}[htbp]
\begin{center}
    \includegraphics[width=6cm]{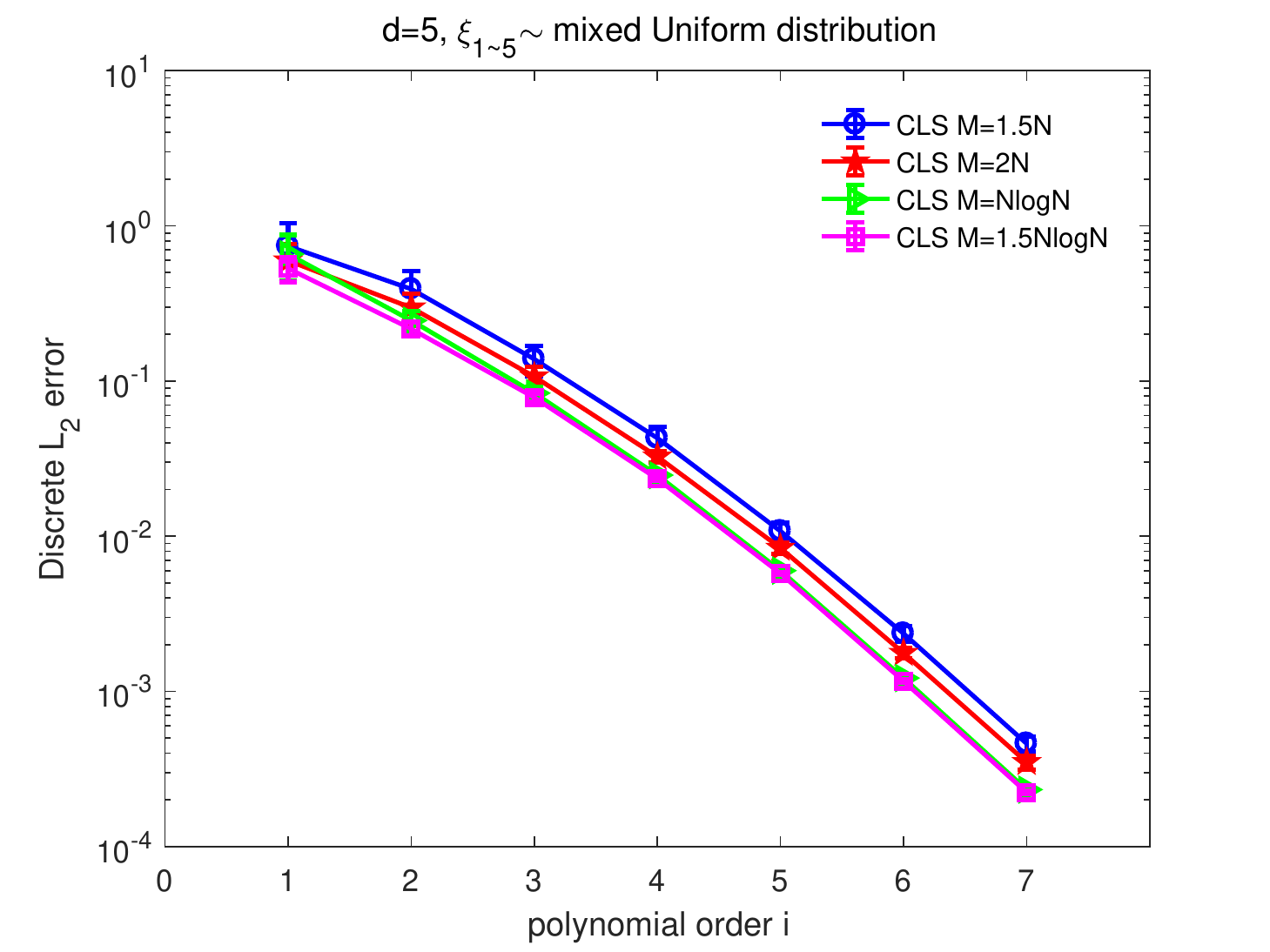}
    \includegraphics[width=6cm]{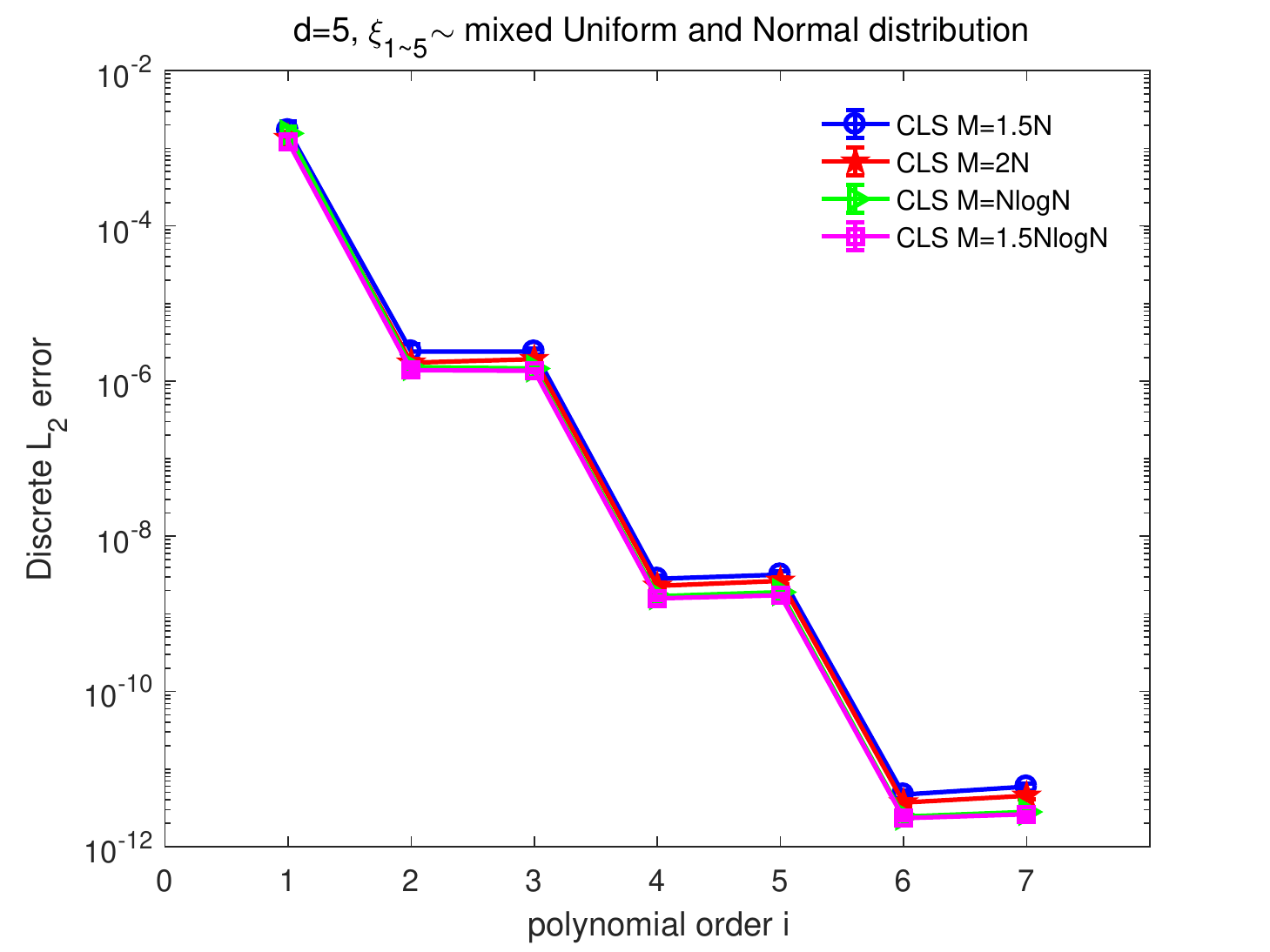}
\end{center}
\caption{Approximation error against polynomial degree for the five dimensional case.}\label{fig:err_d5}
\end{figure}

\begin{figure}[htbp]
\begin{center}
\includegraphics[width=6cm]{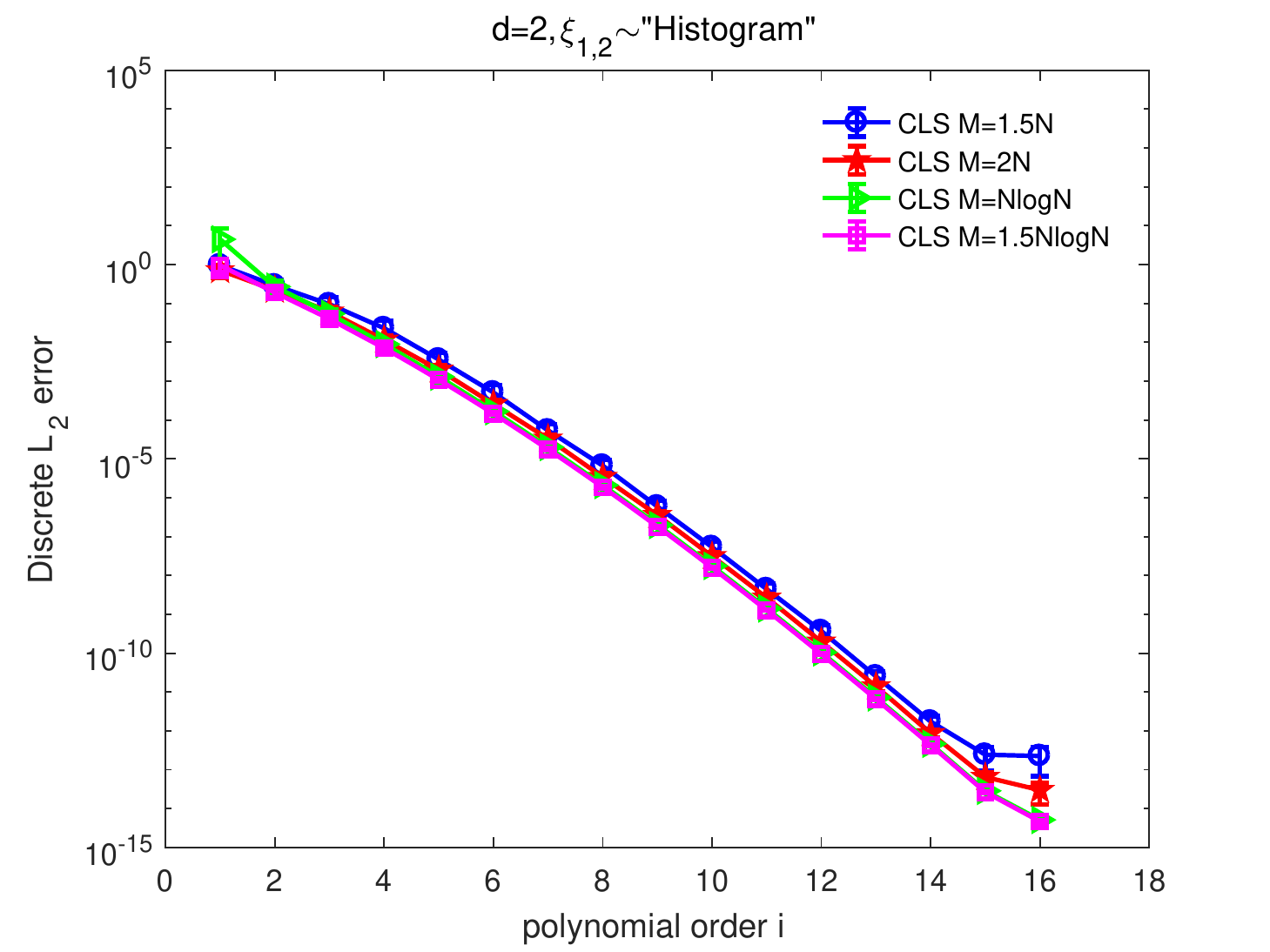}
 \includegraphics[width=6cm]{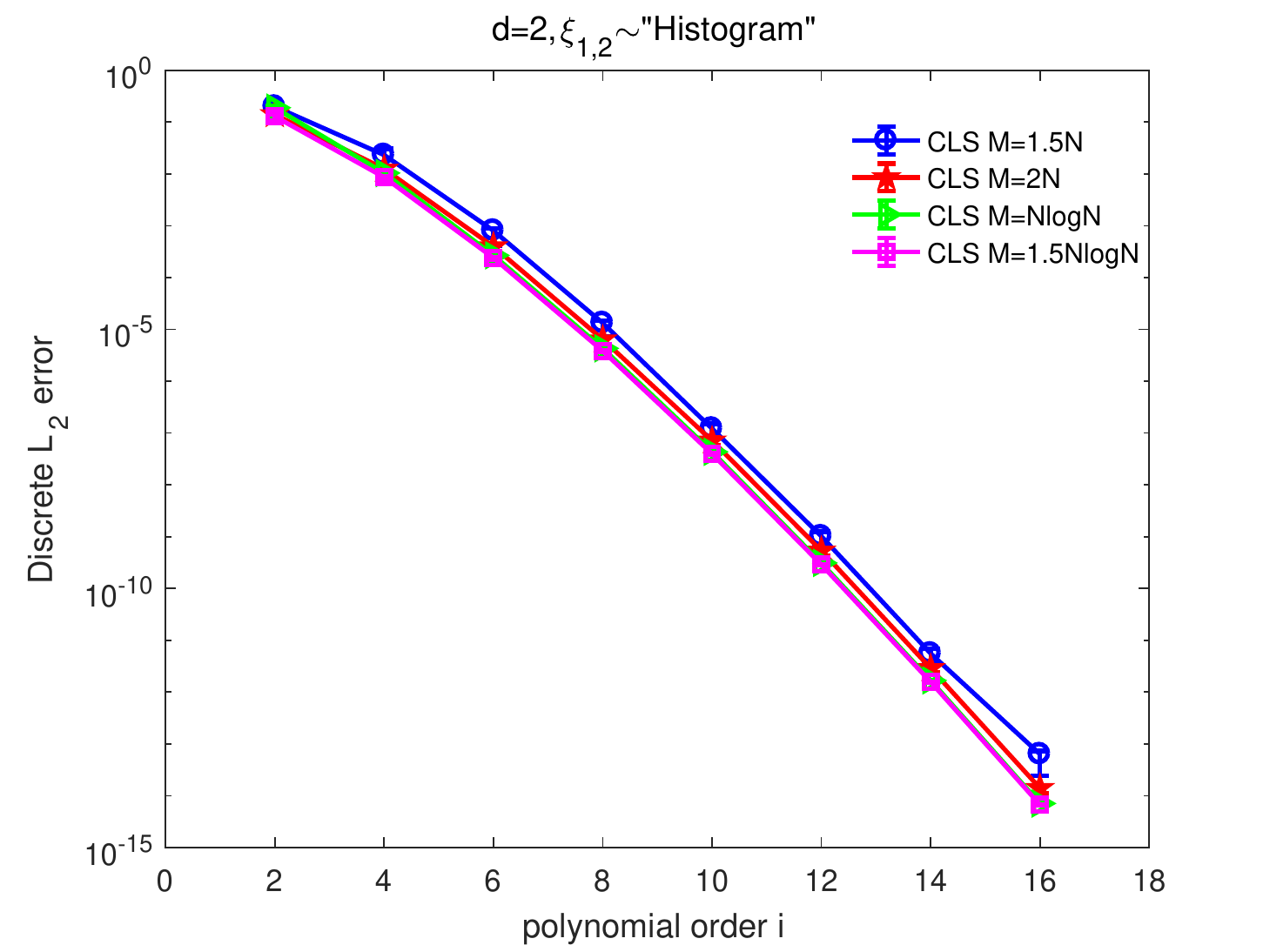}
\\
\includegraphics[width=6cm]{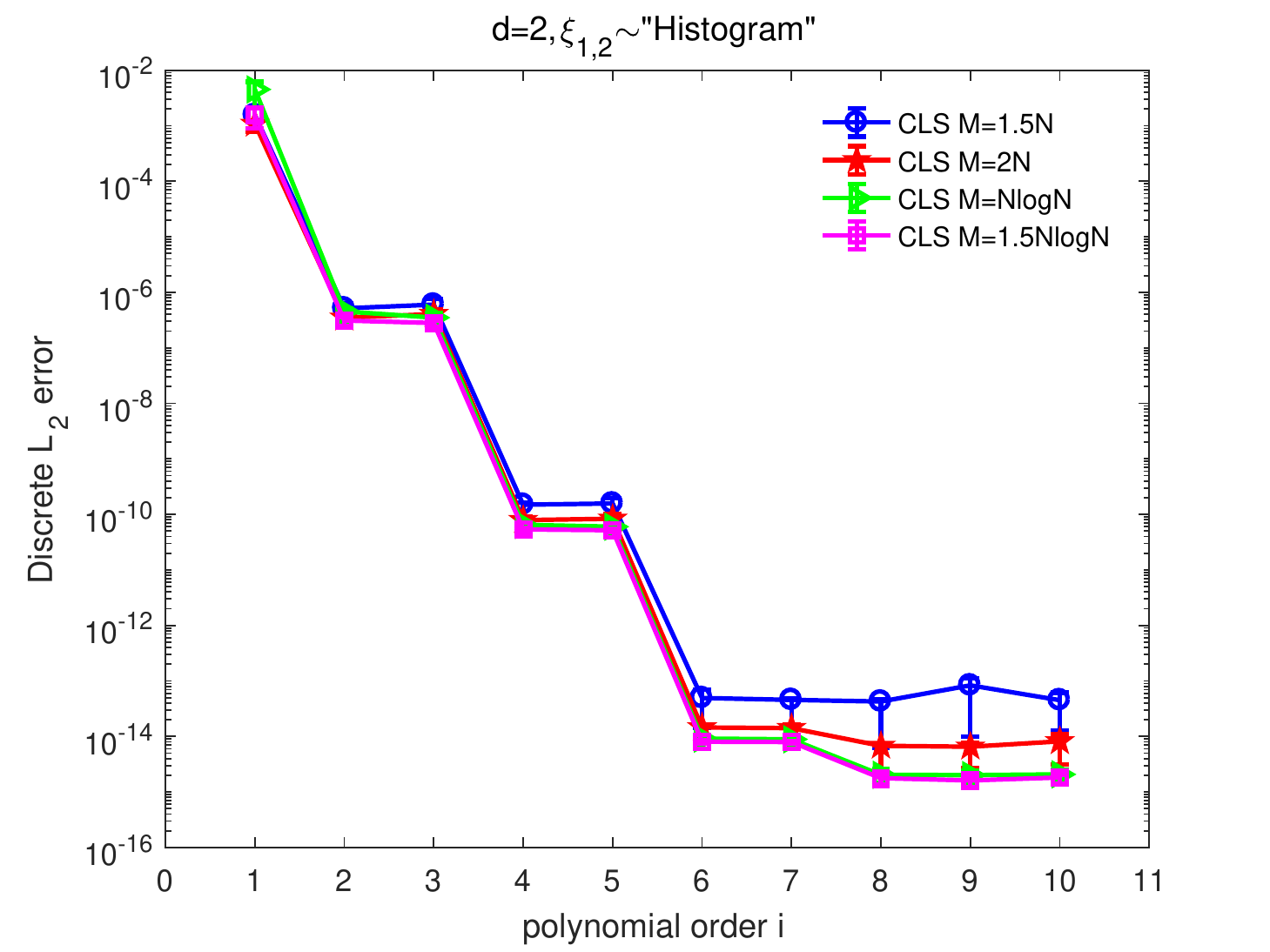}
\includegraphics[width=6cm]{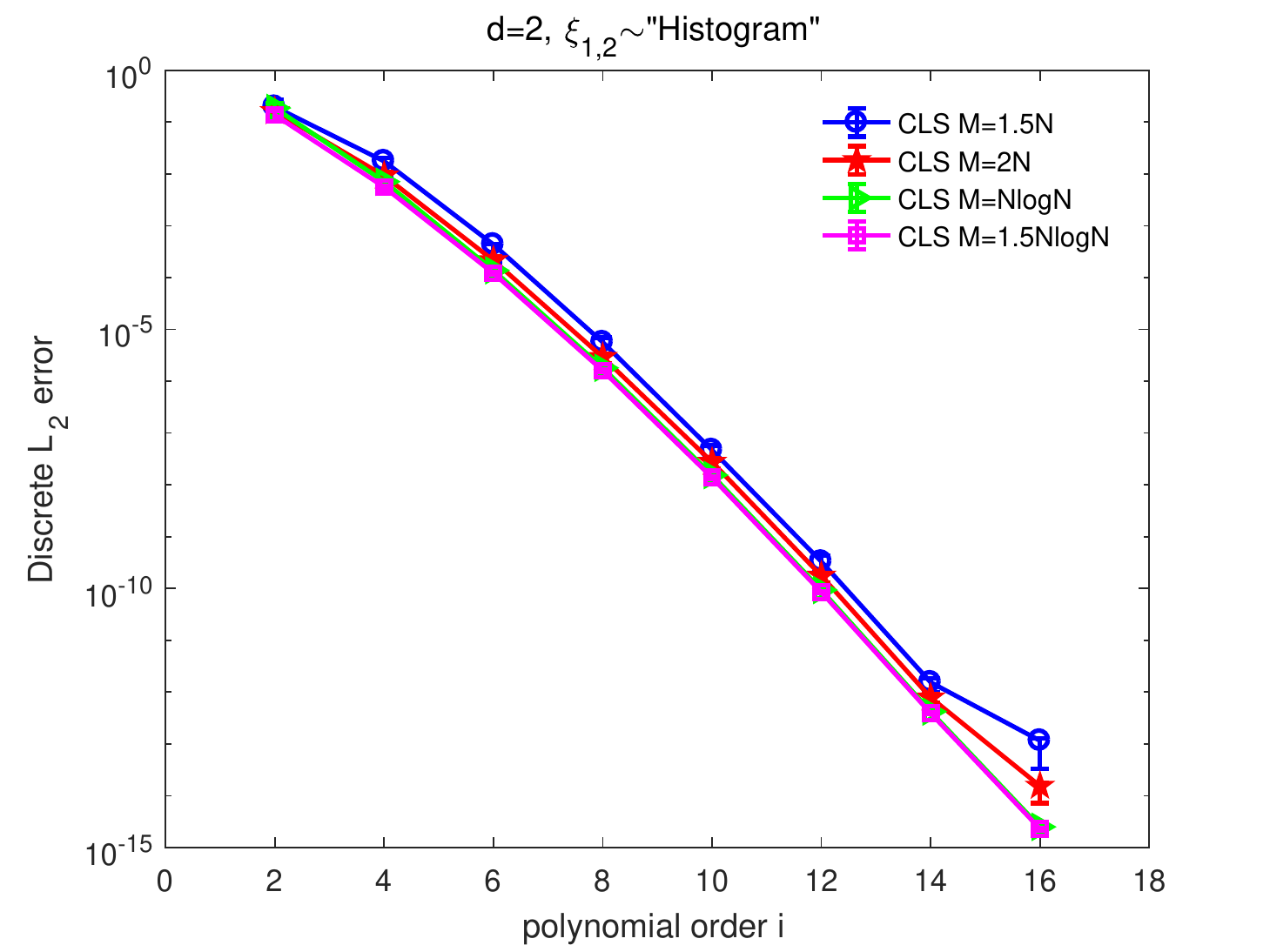}
\end{center}
\caption{Approximation error against polynomial degree for the two dimensional case.}\label{fig:rawdata_ferr_d2}
\end{figure}

\begin{figure}[htbp]
\begin{center}
\includegraphics[width=6cm]{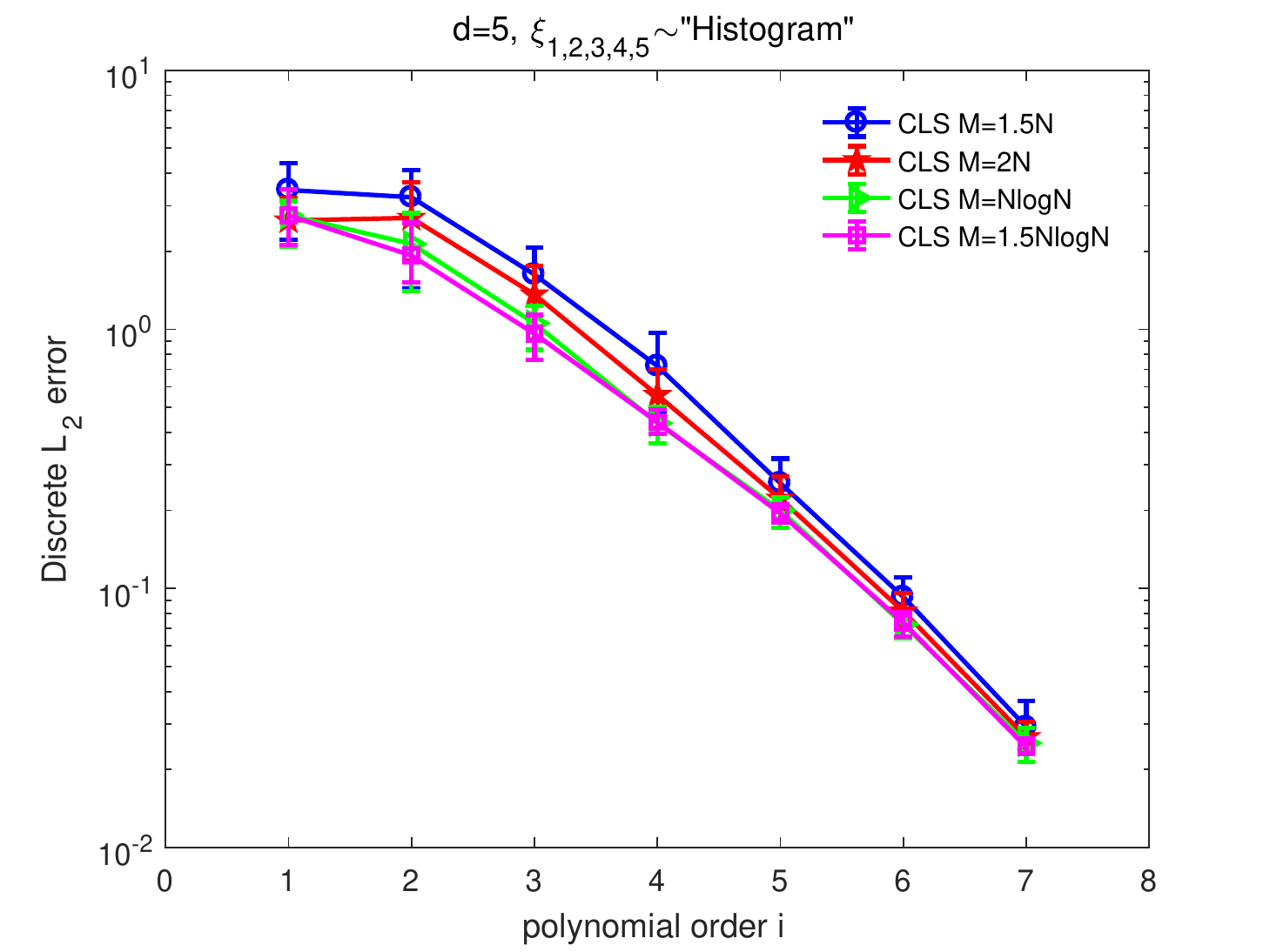}
 \includegraphics[width=6cm]{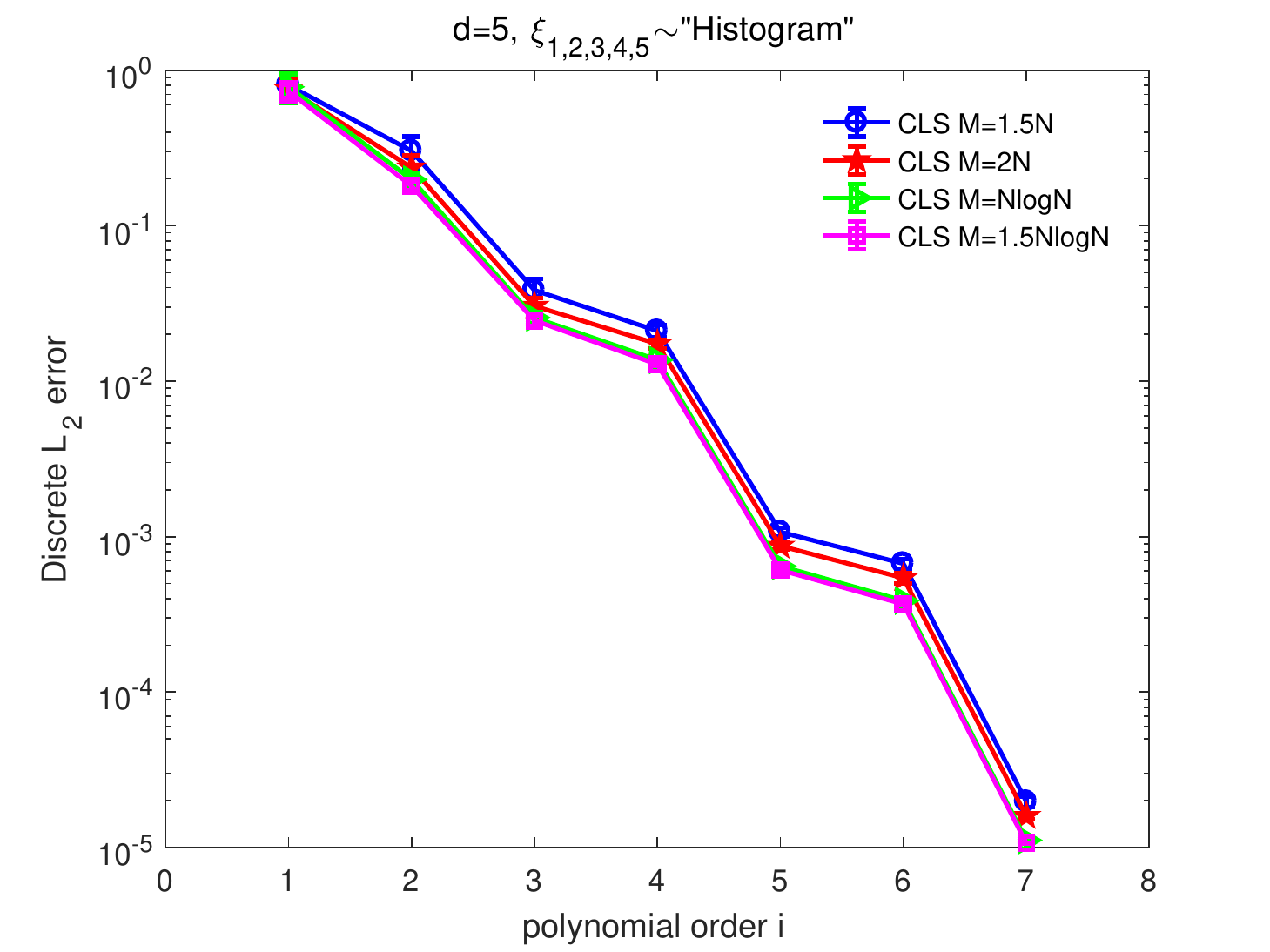}
\\
\includegraphics[width=6cm]{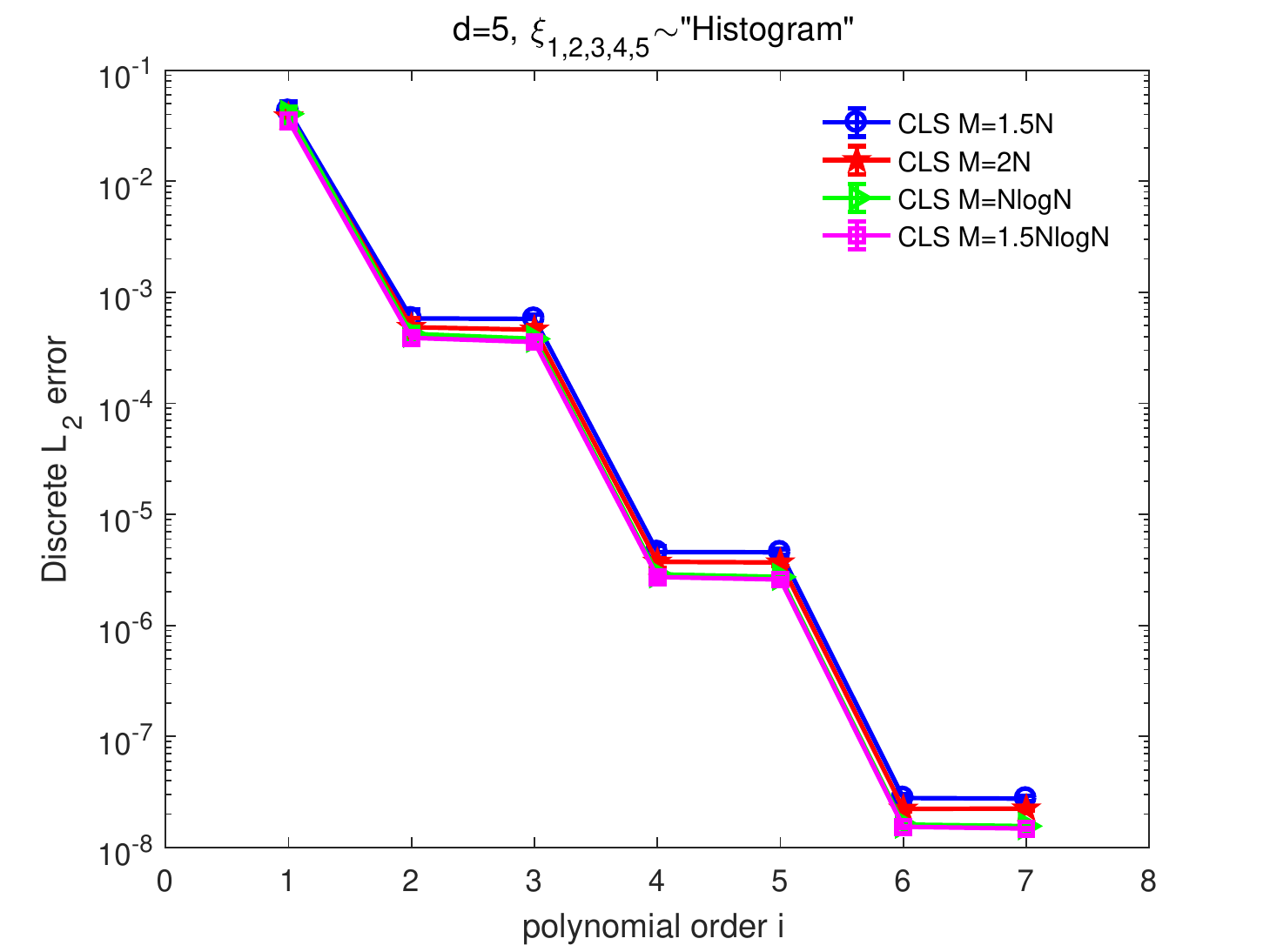}
\includegraphics[width=6cm]{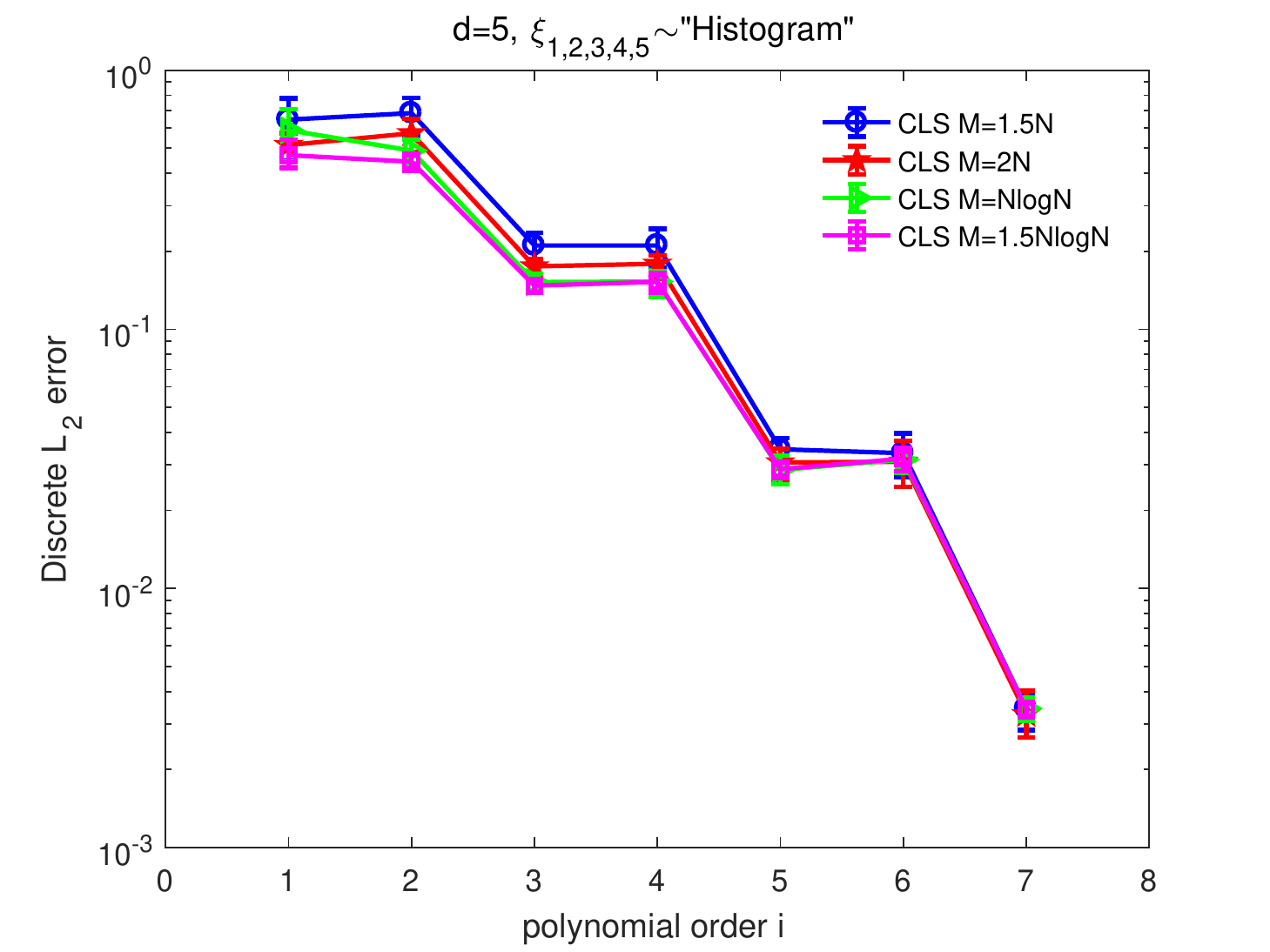}
\end{center}
\caption{Approximation error against polynomial degree for the five dimensional case.}\label{fig:rawdata_ferr_d5}
\end{figure}
\subsubsection{Resistor network}
We now consider a electrical resistor network given in Fig. \ref{fig:resistor_network}. The network is comprised of $d=2p$ resistances $R_i$ of uncertain Ohmage and the network is driven by a voltage source providing a known potential $V_0=1$. We are interested in determining the voltage at $V$, which depends on the $d=2p$ resistances. We set the resistances as random parameters with $d=2$ and $d=4$ cases. To be concrete, we consider the two dimensional parameters uniformly distributed in the interval $\xi_i\in[10,100]$ and the four dimensional parameters with different exponential distribution ($\xi_1\sim \textmd{Exp}(0.9), \xi_2\sim \textmd{Exp}(1.1), \xi_3\sim \textmd{Exp}(0.8), \xi_4\sim \textmd{Exp}(1.0)$). We first use the moments (that are computed with 1000 samples) information to construct a data-driven bases set and then construct the approximation via the weighted least-squares approximation. The accuracy as a function of polynomial order is displayed in Fig. \ref{fig:resistor_network_accuracy}. Similar as in the previous examples, the Christoffel least-squares can provide very stable and accurate approximation results.
\begin{figure}[htbp]
\centering
\includegraphics[width=12cm]{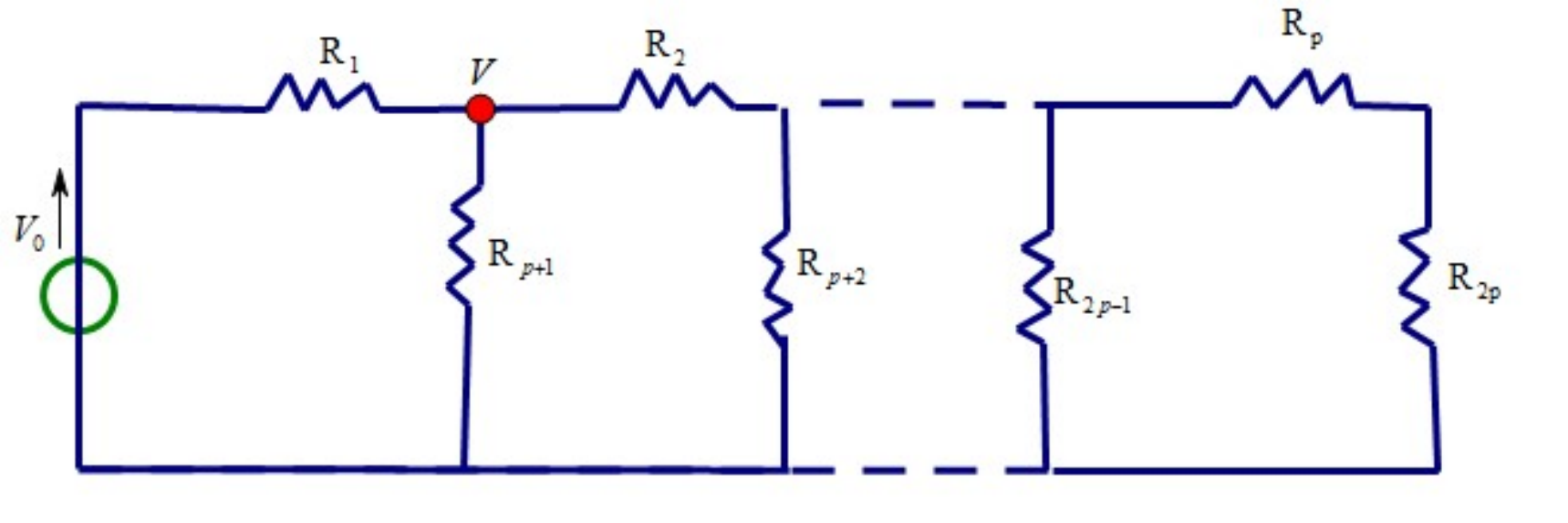}
\caption{Resistor network comprised of $d=2p$ resistances $\{R_i\}_{i=1}^d$ of uncertain ohmage and the network is driven by a voltage source providing a known potential $V_0$.}\label{fig:resistor_network}
\end{figure}

\begin{figure}[htbp]
\begin{center}
\includegraphics[width=6cm]{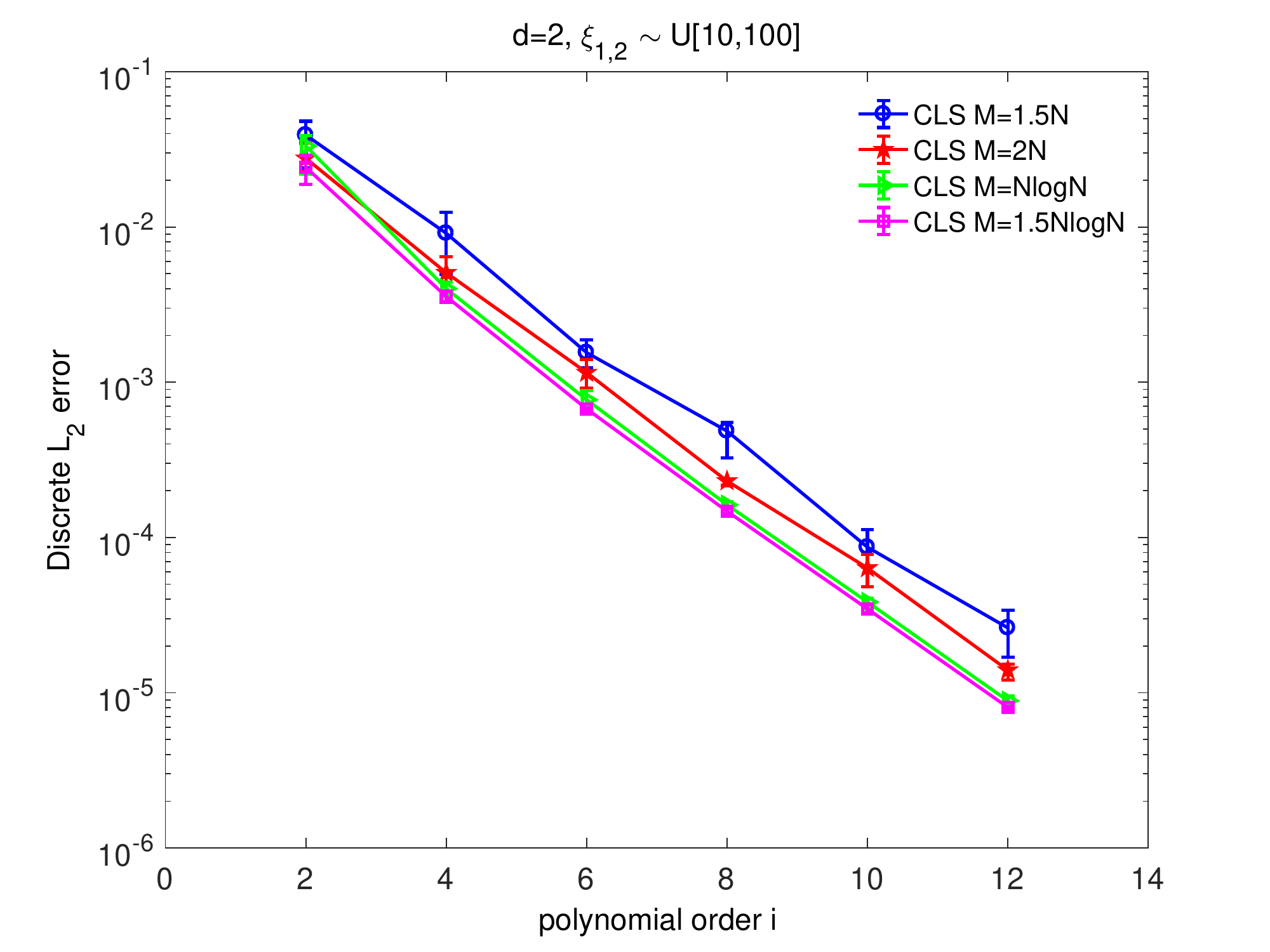}
\includegraphics[width=6cm]{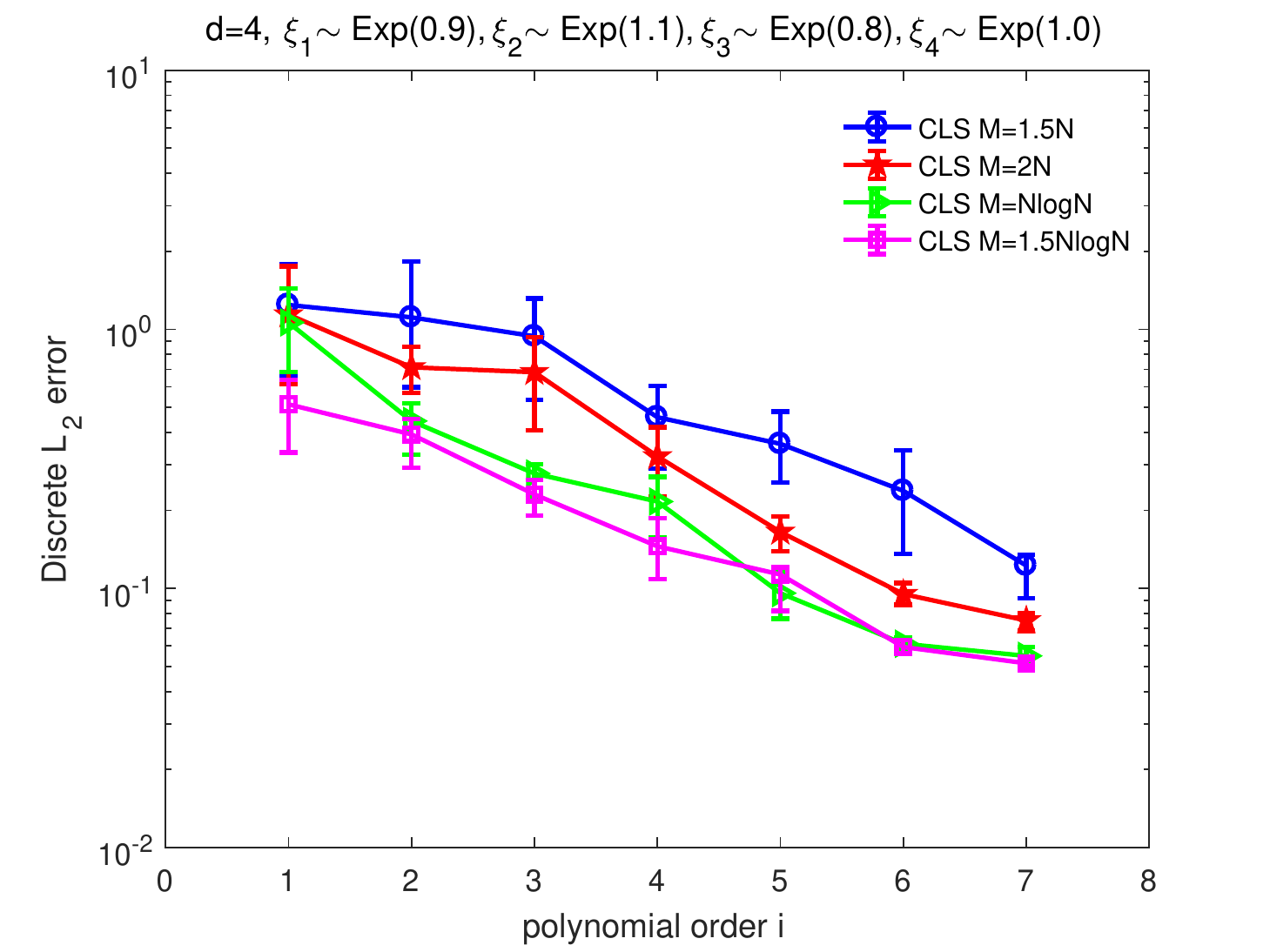}
\end{center}
  \caption{Approximation error against polynomial degree $k$. Left: The two-dimensional isotropic uniform distribution. Right: The four-dimensional anisotropic exponential random distribution.}\label{fig:resistor_network_accuracy}
\end{figure}
\subsubsection{PDEs with random input}
We finnaly consider the following stochastic elliptic PDE
\begin{equation}\label{eq:PDEmodel}
\begin{cases}
-\nabla\cdot(a(\mathbf{y},\omega)\nabla u(\mathbf{y},\omega))=f(\mathbf{y},\omega)\quad \textmd{in} \ \,\mathcal{D}\times\Omega,\\
u(\mathbf{y},\omega)=0 \quad \quad \quad \quad \quad \quad \quad \quad \quad \quad \  \textmd{on} \ \, \partial\mathcal{D}\times\Omega
\end{cases}
\end{equation}
with spatial domain $\mathcal{D}=[0,1]^{2}$. We set a deterministic load $f(\mathbf{y},\omega)=\cos( y_{1})\sin( y_{2})$ for these numerical examples. The random diffusion coefficient $a_{N}(\mathbf{y},\omega)$ is chosen as in  \cite{Babuka_2010scEPDE}:
\begin{equation*}
\log(a_{N}(\mathbf{y},\omega)-0.5)=1+\xi_{1}(\omega)\Big(\frac{\sqrt{\pi}L}{2}\Big)^{1/2}+\sum_{i=2}^{5} \zeta_{i}g_{i}(\mathbf y)\xi_{i}(\omega),
\end{equation*}
where
\begin{equation*}
\zeta_{i}:=(\sqrt{\pi}L)^{1/2}\exp\Big(\frac{-(\lfloor\frac{i}{2}\rfloor\pi L)^{2}}{8}\Big),\, \ \textmd{for } \ \,i>1
\end{equation*}
and
\begin{equation*}
g_{i}(\mathbf{y}):= \begin{cases}
\sin\Big(\frac{-(\lfloor\frac{i}{2}\rfloor\pi  y_{1}}{L_{p}}\Big), \ \ i \ \, \textmd{even},\\[12pt]
\cos\Big(\frac{-(\lfloor\frac{i}{2}\rfloor\pi  y_{1}}{L_{p}}\Big), \ \ i \ \, \textmd{odd}.
\end{cases}
\end{equation*}
Here $\{\xi_{i}\}^{d}_{i=1}$ are independent random variables. For $y_{1}\in[0,1]$, let $L_{c}=1/12$ be a desired physical correlation length for $a(\mathbf{y},\omega)$. Then the parameter $L_{p}$ and $L$ are $L_{p}=\max\{1,2L_{c}\}$ and $L=\frac{L_{c}}{L_{p}}$, respectively. In our numerical test, for each samples, the deterministic elliptic equation are solved by a standard finite element method with a fine mesh. The quantities of interests is the solution $u(y) = u(0.5, 0.5; \xi).$ We set the parametric density as $\xi_1,\xi_2\sim \textmd{Bino}(20,0.5)$ and $\xi_1\sim N(0,1), \xi_2\sim N(0.1,1.2)$. Approximation results are shown in Fig. \ref{fig:PDE_bounded_f54err_d2} with different sampling rates. A good approximation result is observed.
\begin{figure}[htbp]
\begin{center}
    \includegraphics[width=6cm]{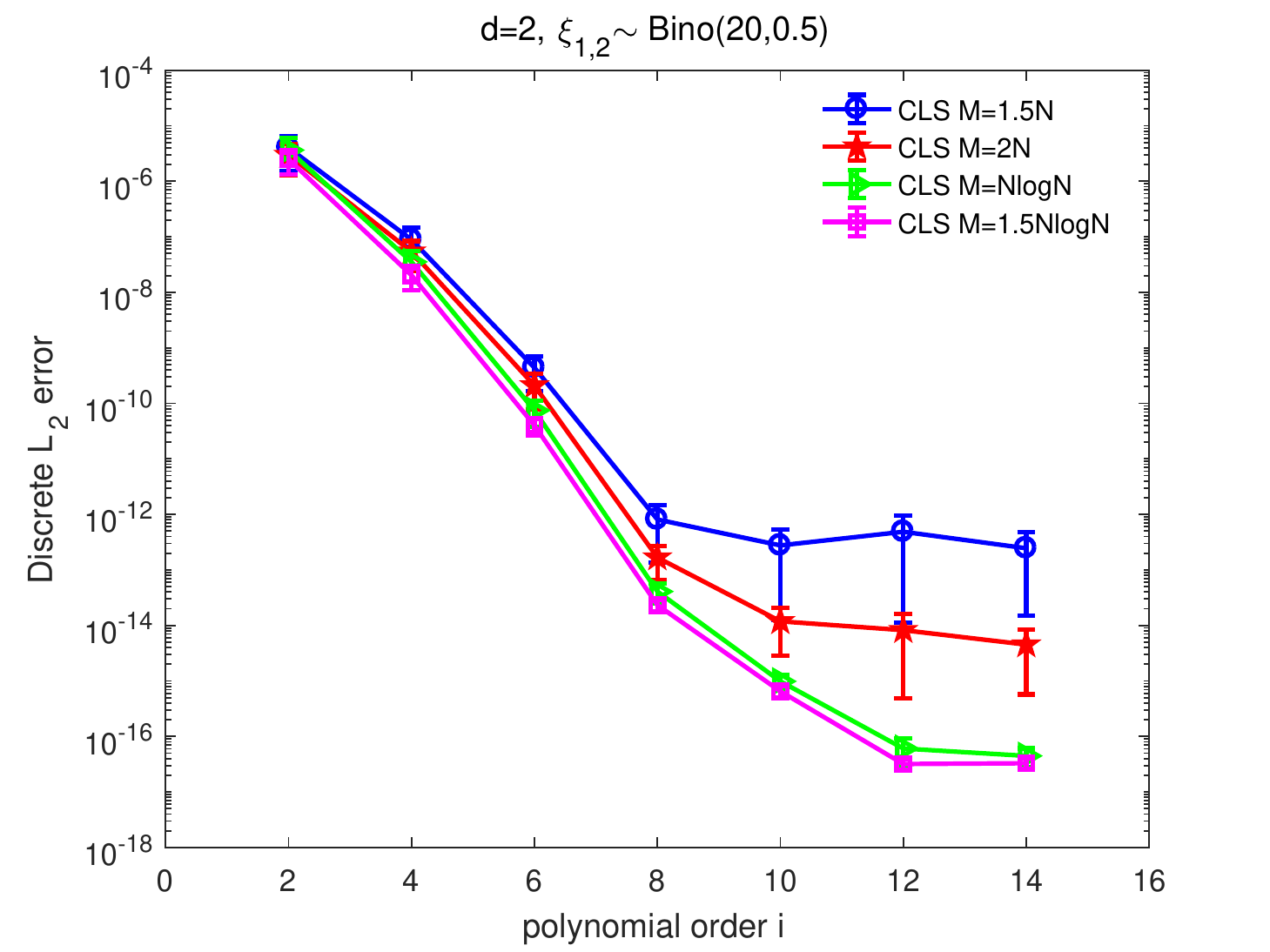}
    \includegraphics[width=6cm]{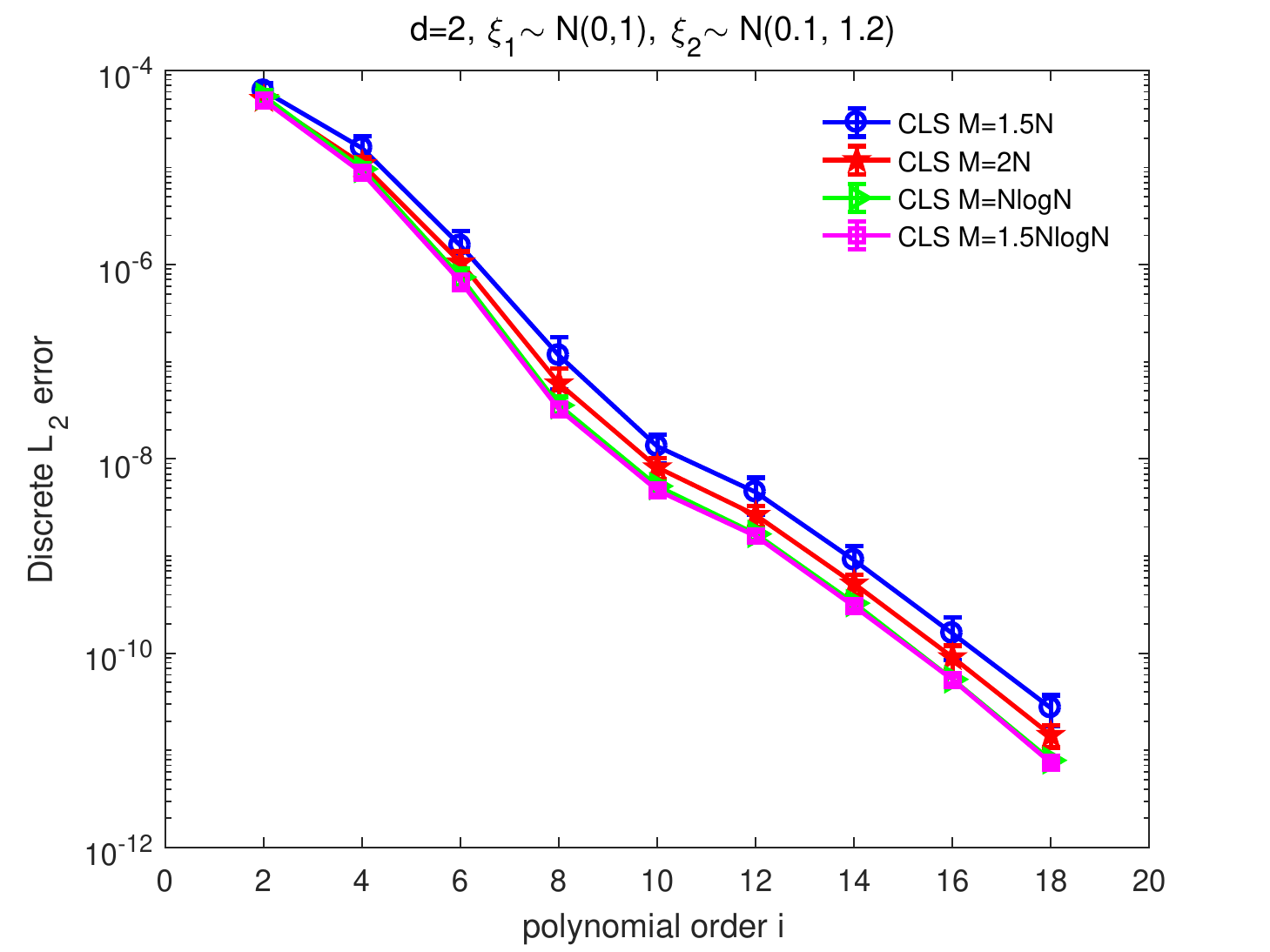}
\end{center}
\caption{Approximation error against polynomial degree of the parametric PDE. }\label{fig:PDE_bounded_f54err_d2}
\end{figure}

\section{Conclusions}
We have combined the idea of  data-driven polynomial chaos expansions with the weighted least-square approach to solve UQ problems. We adopt the bases construction procedure by following \cite{Ahlfeld_2016SAMBA} and then propose to use the weighted least-squares approach to solve UQ problems. Our sampling strategy is independent of the random input. More precisely, we propose to sampling with the equilibrium measure, and this measure is also independent of the data-driven bases. Thus, the procedure can be done in prior (or in a off-line manner). Moreover, the proposed Christoffel function weighted least-squares problem is linearly stable in many cases of interests -- the required number of PDE solvers depends linearly on the number of bases.

There are, however, many unsolved problems related to this topic:
\begin{itemize}
\item Theoretical foundation. As discussed in Section 3.2. The assumption is that the probability density functions are continuous. However, this approach also work well for densities of discrete type as long as their moments exist and the determinant of the moment matrix is strictly positive (see more numerical examples in \cite{Oladyskin_2012Datadriven}). Thus, the relevant theorem for these cases is still open.

\item Density error. We have assumed that only sample locations are given, and all the moments are computed by these finite sample locations, and thus this definitely introduces density error. How to quantify (theoretically) and control this error is of great importance. This is also related to the density sensitivity of the underling model.

\item Unbounded domains. We have provided two simple cases for unbounded domain setting. However, unlike the bounded domain cases, for unbounded cases we need to assume that the type of the density is known (while the associated parameters can be unknown). This is obviously unsatisfactory. Another possible approach to deal with such situations is to truncate the domain into a bounded one (potentially large), and perform the computation in the bounded domain. However, this again introduce the truncation error.

\end{itemize}

We finally close this work by remarking that our strategy can also be used in the compressed sampling setting (or, in the $\ell^1$ approach) \cite{Doostan_2011nonadapted,Jakeman_2016generalizedsample,Guo_2017randomquadratures} and we shall report this in our future studies.

\bibliographystyle{plain}
\bibliography{APCSC_lsqr}

\end{document}